\documentclass[12pt]{article}
\usepackage[final]{epsfig}
\usepackage{graphics}
\usepackage{float}
\usepackage{amsmath}
\usepackage{amsfonts}
\usepackage{latexsym}
\usepackage{amssymb}
\usepackage{graphicx}
\usepackage{url}
\usepackage{epstopdf}
\usepackage{verbatim}
\usepackage[margin=1in]{geometry}
\usepackage{amsthm}

\makeatletter
\newtheorem*{rep@theorem}{\rep@title}
\newcommand{\newreptheorem}[2]{%
\newenvironment{rep#1}[1]{%
 \def\rep@title{#2 \ref{##1}}%
 \begin{rep@theorem}}%
 {\end{rep@theorem}}}
\makeatother

\newtheorem{lemma}{Lemma}[section]
\newtheorem{proposition}[lemma]{Proposition}

\newtheorem{remark}[lemma]{Remark}
\newtheorem{example}[lemma]{Example}
\newtheorem{theorem}[lemma]{Theorem}
\newtheorem{definition}[lemma]{Definition}
\newtheorem{corollary}[lemma]{Corollary}

\newtheorem*{theorem*}{Theorem}
\newreptheorem{theorem}{Theorem}

\newcommand{\proofend}{$\Box$\bigskip}

\newcommand{\N}{{\mathbb N}}
\newcommand{\R}{{\mathbb R}}
\newcommand{\Z}{{\mathbb Z}}

\def\proof{\paragraph{Proof.}}

\DeclareMathOperator*{\argmin}{argmin}
\DeclareMathOperator*{\argmax}{argmax}
\newcommand*{\argminl}{\argmin\limits}
\newcommand*{\argmaxl}{\argmax\limits}

\begin{document}

\title{Fourier--Dedekind sums and an extension of Rademacher reciprocity}

\author{E. Tsukerman\footnote{UC Berkeley, e.tsukerman@berkeley.edu}
}

\date{}
\maketitle
\begin{abstract}
Fourier--Dedekind sums are a generalization of Dedekind sums -- important number-theoretical objects that arise in many areas of mathematics, including lattice point enumeration, signature defects of manifolds and pseudo random number generators. A remarkable feature of Fourier--Dedekind sums is that they satisfy a reciprocity law called Rademacher reciprocity. In this paper, we study several aspects of Fourier--Dedekind sums: properties of general Fourier--Dedekind sums, extensions of the  reciprocity law, average behavior of Fourier--Dedekind sums, and finally, extrema of 2-dimensional Fourier--Dedekind sums.

On properties of general Fourier--Dedekind sums we show that a general Fourier--Dedekind sum is simultaneously a convolution of simpler Fourier--Dedekind sums, and a linear combination of these with integer coefficients. We show that  Fourier--Dedekind sums can be extended naturally to a group under convolution. We introduce ``Reduced Fourier--Dedekind sums", which encapsulate the complexity of a Fourier--Dedekind sum, describe these in terms of generating functions, and give a geometric interpretation.

Next, by finding interrelations among Fourier--Dedekind sums, we extend the range on which Rademacher reciprocity Theorem holds.

We go on to study the average behavior of Fourier--Dedekind sums, showing that the average behavior of a Fourier--Dedekind sum is described concisely by a lower-dimensional, simpler Fourier--Dedekind sum.

Finally, we focus our study on 2-dimensional Fourier--Dedekind sums. We find tight upper and lower bounds on these for a fixed $t$, estimates on the argmax and argmin, and bounds on the sum of their ``reciprocals''. 
\end{abstract}

\section{Introduction}
Dedekind sums are important number-theoretical objects that arise in many areas of mathematics. Among others, these include: geometry (lattice point enumeration in polytopes \cite{BR1}), topology (signature defects of manifolds \cite{HZ1}) and algorithmic complexity (pseudo random number generators \cite{K1}). Fourier--Dedekind sums are generalizations of Dedekind sums that unify many variations of the Dedekind sum that have appeared in the literature. The author's own interest in Fourier--Dedekind sums arose in connection with symplectic embeddings of ellipsoids. It is known that the symplectic invariants called ECH Capacities determine precisely when such an embedding is possible, and these in turn lend themselves to combinatorial analysis. For more information on the topic, see \cite{Hutch1} and \cite{MD1}.

The Fourier--Dedekind sum is defined by
\begin{align} \label{FDDef}
S_{(a_1,a_2,\ldots,a_{d};b)}(n)=\frac{1}{b}\sum_{j=1}^{b-1} \frac{\xi_{b}^{jn}}{(1-\xi_{b}^{ja_1})(1-\xi_{b}^{ja_2})\cdots(1-\xi_{b}^{ja_{d}})},
\end{align}
where $a_1,a_2,\ldots,a_d,b \in \N$, $b>1$ is relatively prime to each $a_i$ and $\xi_b=e^{\frac{2 \pi i}{b}}$. A remarkable feature of Fourier--Dedekind sums is that they satisfy a reciprocity law called Rademacher reciprocity \cite{BDR1}. To state it, we recall the definition of the restricted partition function $p_{\{a_1,a_2,\ldots,a_d\}}(n)$:
\[
p_{\{a_1,a_2,\ldots,a_d\}}(n):=\#\{(m_1,m_2,\ldots,m_d) \in {\Z}^d: \text{all } m_j \geq 0, m_1 a_1+\ldots+m_d a_d=n\}.
\]
Let $\text{poly}_{\{a_1,a_2,\ldots,a_d\}}(n)$ denote the polynomial part of $p_{\{a_1,a_2,\ldots,a_d\}}(n)$, which is explicitly given by the formula \cite{B1}
\[
\text{poly}_{\{a_1,a_2,\ldots,a_d\}}(t)=\frac{1}{a_1 \cdots a_d} \sum_{m=0}^{d-1} \frac{(-1)^m}{(d-1-m)!} \sum_{k_1+\ldots+k_d=m} a_1^{k_1} \cdots a_d^{k_d} \frac{B_{k_1} \cdots B_{k_d}}{k_1! \cdots k_d !}t^{d-1-m},
\]
where $B_j$ denotes the $j$th Bernoulli Number.

 \textit{Rademacher reciprocity} states that given pairwise relatively prime positive integers $a_1,a_2,\ldots,a_d$, for each $n=1,2,\ldots,(a_1+\ldots+a_d-1)$,
\begin{align} \label{RademacherRecip}
S_{(a_2,\ldots,a_{d};a_1)}(n)+S_{(a_1,a_3,a_4,\ldots,a_{d};a_2)}(n)+\ldots+S_{(a_1,a_2,\ldots,a_{d-1};a_d)}(n)=-\text{\textnormal{poly}}_{\{a_1,a_2,\ldots,a_d\}}(-n).
\end{align}
Let
\[
R_{(a_1,a_2,\ldots,a_d)}(t):=S_{(a_2,a_3,\ldots,a_d;a_1)}(t)+S_{(a_1,a_3,a_4,\ldots,a_d;a_2)}(t)+\ldots+S_{(a_1,a_2,\ldots,a_{d-1};a_d)}(t).
\]
The object $R_{(a,1,b)}(t)$ is the main non-trivial ingredient in the Ehrhart quasipolynomial enumerating the lattice points in rational right triangles. As an illustration, we consider right triangles with right angle at the origin. For $e$ and $f$ relatively prime, set 
\[
\mathcal{T}=\{(x,y) \in {\Z}^2 : x \geq 0, y \geq 0, ex+fy \leq r\}.
\]
Then for $t \in {\N}_{\geq 0}$,
\[
L(t):= \#\{(x,y) \in {\Z}^2 : x \geq 0, y \geq 0, ex+fy \leq tr\}
\]
is equal to \cite[Pg. 44]{BR1}
\begin{align} \label{L(t)eq}
L(t)=\frac{1}{2e f}(t r)^2 +\frac{1}{2}(t r)\left( \frac{1}{e}+\frac{1}{f}+\frac{1}{ef}\right)+\frac{1}{4}\left(1+\frac{1}{e}+\frac{1}{f}\right)+\frac{1}{12}\left(\frac{e}{f}+\frac{f}{e}+\frac{1}{ef}\right)+R_{\{e,f,1\}}(-tr).
\end{align}

In particular, the complexity of $L(t)$ lies in $R_{\{e,f,1\}}(-tr)$. 

\section{Main results}

In the first section,``Properties of general Fourier--Dedekind sums", we aim to understand Fourier--Dedekind sums in maximal generality.
Our take-off point is the following theorem, which characterizes Fourier--Dedekind sums in terms of convolutions. Let $T^a$ denote the shift operator by $a$:  $T^a (f(t))=f(t+a)$. The convolution $*$ of two $b$-periodic functions $f$ and $g$ is the function defined by
\[
(f*g)(t)=\sum_{m=0}^{b-1} f(t-m) g(m).
\]
We  denote the indicator function of the integers by $\delta_{\Z}$, so that $\delta_{\Z}(x)=1$ if $x \in \Z$ and $\delta_{\Z}(x)=0$ otherwise.

\begin{reptheorem}{convolutionCharacterization}
Let $f$ and $g$ be $b$-periodic functions with $\sum_{k=0}^{b-1} f(k)=\sum_{k=0}^{b-1} g(k)=0$ (so that $\hat{f}(0)=\hat{g}(0)=0$). The equation 
\begin{align}
(I-T^{a_1})(I-T^{a_2})\cdots (I-T^{a_d})(f * g)=g
\end{align}
holds if and only if 
\begin{align}
\hat{f}(k)=\frac{1}{(1-\xi^{k a_1})(1-\xi^{k a_2})\cdots (1-\xi^{k a_d})}
\end{align}
 whenever $\hat{g}(k) \neq 0$. In particular, if $\hat{g}(k) \neq 0$ for $1\leq k \leq b-1$, then $f(t)=S_{(a_1,a_2,\ldots,a_d;b)}(t)$.
\end{reptheorem}

The next two Theorems show that Fourier--Dedekind sums can be extended naturally into a group under $*$ with simple inverses.

\begin{reptheorem}{functionalEqConv}
The Fourier--Dedekind sum satisfies the functional equation
\[
 (I-T^{a_{d}}) S_{(a_1,a_2,\ldots,a_d;b)}=S_{(a_1,a_2,\ldots,a_{d-1};b)}
\]
and
\[
S_{(a_1,a_2,\ldots,a_d;b)}= S_{(a_1;b)} * S_{(a_2;b)} * \ldots * S_{(a_d;b)}.
\]
\end{reptheorem}

\begin{reptheorem}{group}
The set consisting of all finite products (with operation $*$) of generators of the form 
\[
S_{(a;b)}(t),
\]
\[
 (I-T^a)\delta_{\Z}(\frac{t}{b}),
\]
and
\[ 
  S_b(t),
  \]
 with $(a,b)=1$, is an abelian group with identity $S_b$. The inverse of $S_{(a;b)}$ is $(I-T^a)\delta_{\Z}$. 
\end{reptheorem}

As a consequence, we are able to show that a Fourier--Dedekind sum is a $\Z$-linear combination of simpler Fourier--Dedekind sums:

\begin{reptheorem}{linearComb} Let $d \geq 1$. Then
\begin{align}
b S_{(a_1,a_2,\ldots,a_d)}(t)=- \sum_{k=1}^{b-1} k S_{(a_1,a_2,\ldots,a_{d-1})}(t+k a_d)
\end{align}
\end{reptheorem}

We introduce the Reduced Fourier--Dedekind sum $\tilde{S}_{(a_1,a_2,\ldots,a_d;b)}(t)$, defined by
\[
\tilde{S}_{(a_1,a_2,\ldots,a_d;b)}(t):=\sum_{\substack{1 \leq k_1, k_2,\ldots, k_d \leq b-1,\\ a_1 k_1+a_2 k_2 +\ldots+a_d k_d \equiv -t \pmod{b}}} k_1 k_2 \cdots k_d.
\]

The Reduced Fourier--Dedekind sum is the nontrivial part of a Fourier--Dedekind sum:

\[
 S_{(a_1,a_2,\ldots,a_d;b)}(t)=
\frac{(-1)^d}{b^d}[\tilde{S}_{(a_1,a_2,\ldots,a_d;b)}(t)-\frac{1}{b}\binom{b}{2}^d].
\]

We go on to describe Reduced Fourier--Dedekind sums in  terms of simple generating functions.

\begin{reptheorem}{genFuncThm}
For any $a_1,a_2,\ldots,a_d,b \in \N$, \\\\
$
 \quad \quad \quad \tilde{S}_{(a_1,a_2,\ldots,a_d;b)}(t)=
$
\[
 \sum_{j= -\infty}^{\infty} [z^{-t+bj}] \left(z^{a_1}+2z^{2a_1}+\ldots+(b-1)z^{(b-1)a_1}\right)\cdots \left(z^{a_d}+2z^{2a_d}+\ldots+(b-1)z^{(b-1)a_d}\right ).
\]
\end{reptheorem}

This allows us to interpret Reduced Fourier--Dedekind sums geometrically. More specifically, we consider a torus $T$ whose fundamental domain is $F=\{(x_1,x_2,\ldots,x_n) \in \R^n: 0 \leq x_1,x_2,\ldots,x_n < b\}$. We assign each point $(x_1,x_2,\ldots,x_n)$ in $F$ the weight $x_1  x_2 \cdots x_n$, which one could interpret as a suitable volume, and extend periodically to $T$. We let $H$ be the hyperplane $a_1x_1+a_2x_2+\ldots+a_nx_n \equiv -t \pmod{b}$. The value  $\tilde{S}_{(a_1,a_2,\ldots,a_d;b)}(t)$ is equal to the weighted sum over lattice points in $T \cap H$ (see Figure \ref{squareLattice} for the 2-dimensional case).

In the second section, titled ``An extension of Rademacher reciprocity", we extend the range of possible values $n$ under which the Theorem of Rademacher reciprocity holds.

\begin{reptheorem}{RademacherExtended}Let $a_1,a_2,\ldots,a_d \in \N$ be pairwise relatively prime. Let $n \in \Z$. If one of $(i)-(iii)$ holds, where
\begin{itemize}
\item[(i).] $1-\min\{a_1,a_2,\ldots,a_d\} \leq n \leq -1$
\item[(ii).] $1 \leq n \leq a_1+a_2+\ldots+a_d-1$
\item[(iii).] $a_1+a_2+\ldots+a_d+1 \leq n \leq a_1+a_2+\ldots+a_d+\min\{a_1,a_2,\ldots,a_d\}-1$
\end{itemize}
then
\[
S_{(a_2,\ldots,a_{d};a_1)}(n)+S_{(a_1,a_3,a_4,\ldots,a_{d};a_2)}(n)+\ldots+S_{(a_1,a_2,\ldots,a_{d-1};a_d)}(n)=-\text{\textnormal{poly}}_{\{a_1,a_2,\ldots,a_d\}}(-n).
\]
\end{reptheorem}

In the third section, ``Average behavior of Fourier--Dedekind sums", we study the average behavior of a Fourier--Dedekind sum as the $a_i$'s vary.  

More precisely, we define the average over the $i$th variable of $S_{(a_1,a_2,\ldots,a_{d-1},a_d;b)}$ at $t$, denoted by $S_{(a_1,a_2,\ldots,\bar{a}_i,\ldots,a_d;b)}(t)$,  to be
\[
S_{(a_1,a_2,\ldots,\bar{a}_i,\ldots,a_d;b)}(t):=\frac{1}{\phi(b)} \sum_{\substack{1 \leq m \leq b-1\\ (m,b)=1}}  S_{(a_1,a_2,\ldots,a_{i-1},m,a_{i+1},\ldots,a_d;b)}(t).
\] 
The average over all variables of $S_{(a_1,a_2,\ldots,a_{d-1},a_d;b)}$ at $t$, denoted by $S_{(\bar{a}_1,\bar{a}_2,\ldots,\bar{a}_d;b)}(t)$, is defined to be
\[
S_{(\bar{a}_1,\bar{a}_2,\ldots,\bar{a}_d;b)}(t):=(\frac{1}{\phi(b)})^d \sum_{\substack{1 \leq m_1,m_2,\ldots,m_d \leq b-1\\ (m_i,b)=1}}  S_{(m_1,m_2,\ldots,m_d;b)}(t)
\]

\begin{reptheorem}{avg} Let $b \geq 3$ and let $(a_i,b)=1$ for each $i$. For every $t \in \Z$,
\[
S_{(a_1,a_2,\ldots,a_{d-1},\bar{a}_d;b)}(t)=\frac{1}{2} S_{(a_1,a_2,\ldots,a_{d-1})}(t).
\] 
and
\[
S_{(\bar{a}_1,\bar{a}_2,\ldots,\bar{a}_d;b)}(t)=\frac{1}{2^d} S_{b}(t)=\frac{\delta_{\Z}(\frac{t}{b})-\frac{1}{b}}{2^d} .
\] 
\end{reptheorem}

In the final section, ``Bounds, maxima and minima of 2-dimensional Fourier--Dedekind sums" we focus on $2$-dimensional Fourier--Dedekind sums. Our three main results in this section are the following:

\begin{reptheorem}{DedekindSumBounds} For all $a_1,a_2$ coprime to $b$,

\[
-\frac{(b-1)(b-5)}{12b} \leq S_{(a_1,a_2;b)}(0) \leq \frac{(b-1)(b+1)}{12 b}.
\]
The upper bound holds if and only if $a_1+a_2 \equiv 0 \pmod{b}$. The lower bound holds if and only if $a_1 \equiv a_2 \pmod{b}$.

For all $a_1,a_2$ coprime to $b$ and $1 \leq t \leq b-1$,
\[
-\frac{(b-1)(b+1)}{12b} \leq S_{(a_1,a_2;b)}(t) \leq \frac{(b-1)(b-5)}{12 b}.
\]
The upper bound holds if and only if $a_1 \equiv -a_2 \equiv t \pmod{b}$. The lower bound holds if and only if $a_1 \equiv a_2 \equiv t \pmod{b}$.

\end{reptheorem}

To understand the location of maxima and minima, it suffices to assume that $a_2=1$ by a change of variable.

\begin{reptheorem}{Optimization} For $a,b,t \in \N$,
\[
\argmaxl_{1 \leq t \leq b} S_{(a,1;b)}(t) \subset [\frac{b+1}{2},\frac{b+1}{2}+a].
\]
and
\[
\argminl_{1 \leq t \leq b} S_{(a,1;b)}(t) \subset [1,\min\{a,\frac{b+1}{2}\}].
\]
\end{reptheorem}

Finally, we bound the reciprocal sum. The following result yields polynomial bounds on $R_{\{a,1,b\}}(t+a+b)$ whenever polynomial bounds are known for $R_{\{a,1,b\}}(t)$ (e.g., when Rademacher reciprocity holds).

\begin{reptheorem}{boundstaplusb} For every $t \in \Z$,
\[
|R_{\{a,1,b\}}(t+a+b)-R_{\{a,1,b\}}(t)| \leq 1-\frac{1}{2}(\frac{1}{a}+\frac{1}{b}).
\]
\end{reptheorem}

\section{Properties of general Fourier--Dedekind sums}

Throughout this section, unless otherwise stated, we assume that $a_1,a_2,\ldots,a_d \in \N$ are coprime to $b$. Indeed, otherwise, definition (\ref{FDDef}) will include a summand with zero in the denominator. We let $\delta_{\Z}$ denote the indicator function of the integers, so that $\delta_{\Z}(x)=1$ if $x \in \Z$ and $\delta_{\Z}(x)=0$ otherwise.

In preparation for our first result, we recall some basic facts from Fourier Analysis of periodic functions on $\Z$. The convolution $*$ of two $b$-periodic functions $f$ and $g$ is the function defined by
\[
(f*g)(t)=\sum_{m=0}^{b-1} f(t-m) g(m).
\]
Let $\xi=e^{\frac{2 \pi i}{b}}$. Any periodic function $a(n)$ on $\Z$ with period $b$ has a unique discrete Fourier expansion,
\[
a(n)=\sum_{k=0}^{b-1} \hat{a}(k) \xi^{nk},
\]
where
\[
\hat{a}(n)=\frac{1}{b} \sum_{k=0}^{b-1} a(k) \xi^{-nk}
\]
are the Fourier coefficients.

Let $T^a$ denote the shift operator by $a$: $T^a (f(t))=f(t+a)$. The following theorem characterizes Fourier--Dedekind sums in terms of a very general convolution equation.

\begin{theorem} \label{convolutionCharacterization}
Let $f$ and $g$ be $b$-periodic functions with $\sum_{k=0}^{b-1} f(k)=\sum_{k=0}^{b-1} g(k)=0$ (so that $\hat{f}(0)=\hat{g}(0)=0$). The equation 
\begin{align} \label{UPeq}
(I-T^{a_1})(I-T^{a_2})\cdots (I-T^{a_d})(f * g)=g
\end{align}
holds if and only if 
\begin{align} \label{UPeq2}
\hat{f}(k)=\frac{1}{(1-\xi^{k a_1})(1-\xi^{k a_2})\cdots (1-\xi^{k a_d})}
\end{align}
 whenever $\hat{g}(k) \neq 0$. In particular, if $\hat{g}(k) \neq 0$ for $1\leq k \leq b-1$, then $f(t)=S_{(a_1,a_2,\ldots,a_d;b)}(t)$.
\end{theorem}

The assumption $\sum_{k=0}^{b-1} f(k)=\sum_{k=0}^{b-1} g(k)=0$ serves to normalize $f$ and $g$. To prove this result, we will employ the Convolution Theorem for finite Fourier series, which states that 
\[
(f * g)(t) =\frac{1}{b} \sum_{k=0}^{b-1} \hat{f}(k) \hat{g}(k) \xi^{k t}.
\]

\proof
Assume first that the equation holds. Let $c_j \in \Z$, $j=0,1,\ldots,a_1+a_2+\ldots+a_d$, satisfy $\sum_{j=0}^{a_1+a_2+\ldots+a_d} c_j T^j=\prod_{j=1}^{d} (I-T^{a_j})$. Then
\[
(I-T^{a_1})(I-T^{a_2})\cdots (I-T^{a_d})(f * g)(t)=\sum_{j=0}^{a_1+a_2+\ldots+a_d} c_j (f * g)(t+j).
\]
We apply the discrete Fourier transform to both sides of equation (\ref{UPeq}) and use the Convolution Theorem:
\[
\sum_{j=0}^{a_1+a_2+\ldots+a_d} c_j (\frac{1}{b} \sum_{k=0}^{b-1} \hat{f}(k) \hat{g}(k) \xi^{k t} \xi^{k j})  =\frac{1}{b} \sum_{k=0}^{b-1} \hat{g}(k) \xi^{k t}.
\]
Multiplying both sides by $b$ and rearranging,
\[
\sum_{k=0}^{b-1} \hat{g}(k) \xi^{k t}=\sum_{k=0}^{b-1}  \hat{f}(k) \hat{g}(k) \xi^{k t} \sum_{j=0}^{a_1+a_2+\ldots+a_d} c_j \xi^{k j} .
\]
\[
=\sum_{k=0}^{b-1}  \hat{f}(k) \hat{g}(k) \xi^{k t} [\prod_{j=1}^{d} (1-\xi^{a_j k}) ] = \sum_{k=0}^{b-1}   [\prod_{j=1}^{d} (1-\xi^{a_j k}) ] \hat{f}(k) \hat{g}(k) \xi^{k t}.
\]

By uniqueness of the discrete Fourier transform, for each $k \in \Z$,
\[
[\prod_{j=1}^{d} (1-\xi^{a_j k}) ] \hat{f}(k) \hat{g}(k) =  \hat{g}(k).
\]
For the converse, we observe that (\ref{UPeq2}) implies
\[
\frac{1}{b} \sum_{k=0}^{b-1}   [\prod_{j=1}^{d} (1-\xi^{a_j k}) ] \hat{f}(k) \hat{g}(k) \xi^{k t}=\frac{1}{b} \sum_{k=0}^{b-1}\hat{g}(k) \xi^{k t}.
\]
Taking the inverse discrete Fourier transform and applying the Convolution Theorem completes the proof.
\proofend

\begin{definition}
We call the number $d$ in (\ref{FDDef}) the \textit{dimension} of the Fourier--Dedekind sum. 
\end{definition}
In this terminology, the following corollary shows that a Fourier--dedekind sum may be built-up by convolving the more elementary lower-dimensional Fourier--Dedekind sums and  that the operators $I-T^{a_i}$, which yield suitable finite differences, lower the dimension of a Fourier-Dedekind sum. Later we will show that the operation may be reversed in the sense that a $d+1$-dimensional Fourier--Dedekind sum is a $\Z$-linear combination of $d$-dimensional Fourier--Dedekind sums.

\begin{theorem} \label{functionalEqConv}
The Fourier--Dedekind sum satisfies the functional equation
\[
 (I-T^{a_{d}}) S_{(a_1,a_2,\ldots,a_d;b)}=S_{(a_1,a_2,\ldots,a_{d-1};b)}
\]
and
\[
S_{(a_1,a_2,\ldots,a_d;b)}= S_{(a_1;b)} * S_{(a_2;b)} * \ldots * S_{(a_d;b)}.
\]
\end{theorem}

\proof
Both statements follows from the uniqueness of Theorem \ref{convolutionCharacterization} .
\proofend

As an interesting consequence, we note that  for any $i,j \in \{1,2,\ldots,d\}$,
\[
\sum_{k=0}^{\frac{\text{lcm} (a_i,a_j)}{a_i}-1} S_{(a_1,a_2,\ldots,\hat{a}_i,\ldots, a_d;b)}(t+k a_i)=\sum_{k=0}^{\frac{\text{lcm} (a_i,a_j)}{a_j}-1} S_{(a_1,a_2,\ldots,\hat{a}_j,\ldots, a_d;b)}(t+k a_j).
\]

Indeed, by Theorem \ref{functionalEqConv}, we have 
\[
\sum_{k=0}^{\frac{\text{lcm} (a_i,a_j)}{a_i}-1} S_{(a_1,a_2,\ldots,\hat{a}_i,\ldots, a_d;b)}(t+k a_i)=S_{(a_1,a_2,\ldots, a_d;b)}(t)-S_{(a_1,a_2,\ldots, a_d;b)}(t+\text{lcm} (a_i,a_j))
\]
\[
=\sum_{k=0}^{\frac{\text{lcm} (a_i,a_j)}{a_j}-1} S_{(a_1,a_2,\ldots,\hat{a}_j,\ldots, a_d;b)}(t+k a_j).
\]
\proofend

For future reference, we record the following easy result, which concerns the $0$-dimensional Fourier--Dedekind sum.
\begin{lemma} \label{geometricSeries} Let $b,t \in \N$.
\[
S_b(t):=\frac{1}{b} \sum_{j=1}^{b-1} \xi_b^{j t}= \begin{cases} 1-\frac{1}{b} &\mbox{if } t \equiv 0 \pmod{b} \\
-\frac{1}{b} & \mbox{if } t \not \equiv 0 \pmod{b} \end{cases}=\delta_{\Z}(\frac{t}{b})-\frac{1}{b}.
\]
\end{lemma}

\proof
If $t \equiv 0 \pmod{b}$, then the statement is obvious. If $ t \not \equiv 0 \pmod{b}$ then the expression is a geometric series.
\proofend

Our aim now is to extend Fourier--Dedekind sums into a group under the operation $*$. 
Let $V_b$ denote the vector space of real-valued $b$-periodic functions on $\Z$. We define the subspace 
\[
V_b^0:=\{f \in V_b : \frac{1}{b} \sum_{j=0}^{b-1} f(j)=\hat{f}(0)=0\}.
\]
Clearly, $S_{(a_1,a_2,\ldots,a_d;b)}(t) \in V_b^0$. 

\begin{lemma} \label{inverseLemma}
The Fourier--Dedekind sum $S_b$ is the identity in $(V_b^0,*)$. Any Fourier--Dedekind sum has a unique inverse under $*$ in $V_b^0$.
When $d \geq 1$, the inverse of $S_{(a_1,a_2,\ldots,a_d;b)}(t)$  is given by 
\begin{align} \label{inverse}
(I-T^{a_1})(I-T^{a_2})\cdots(I-T^{a_d})\delta_{\Z}(\frac{t}{b}).
\end{align}
\end{lemma}

For the proof, we will employ the following fact. For $k \in \Z$,
\[
T^{k} (f * g)=(T^{k} f) * g=f * (T^{k} g).
\]

\proof
The first statement is straightforward, and the second follows from the Convolution Theorem using the observation that all Fourier coefficients except the first are nonzero.
To see that expression (\ref{inverse}) is an inverse, we employ Theorem \ref{functionalEqConv}:
\[
S_{(a_1,a_2,\ldots,a_d;b)}*(I-T^{a_1})(I-T^{a_2})\cdots(I-T^{a_d})\delta_{\Z}(\frac{t}{b})
\]
\[
=(I-T^{a_1})(I-T^{a_2})\cdots(I-T^{a_d})S_{(a_1,a_2,\ldots,a_d;b)} *  \delta_{\Z}(\frac{t}{b})
\]
\[
=S_b * \delta_{\Z}(\frac{t}{b})= S_b.
\]
\proofend

To summarize,

\begin{theorem} \label{group}
The set consisting of all finite products (with operation $*$) of generators of the form 
\[
S_{(a;b)}(t),
\]
\[
 (I-T^a)\delta_{\Z}(\frac{t}{b}),
\]
and
\[ 
  S_b(t),
  \]
 with $(a,b)=1$, is an abelian group with identity $S_b$. The inverse of $S_{(a;b)}$ is $(I-T^a)\delta_{\Z}$. 
\end{theorem}

 We denote this group by $\mathcal{FD}$. Since we may reduce $a$ modulo $b$, this group is finitely generated. In particular,
\[
\mathcal{FD} \simeq {\Z}^{\phi(b)}.
\] 
In the future, we will employ the notation 
\[
S^{-1}_{(a_1,a_2,\ldots,a_d;b)}(t):=(I-T^{a_1})(I-T^{a_2})\cdots (I-T^{a_d})\delta_{\Z}(\frac{t}{b}).
\]

\begin{corollary} \label{systemOfEquations}
The vector 
\[
(S_{(a_1,a_2,\ldots,a_d;b)}(0),S_{(a_1,a_2,\ldots,a_d;b)}(1),\ldots,S_{(a_1,a_2,\ldots,a_d;b)}(b-1))^t
\]
is the unique solution to the system of equations
\begin{align} \label{matrix} \left( \begin{array}{cccc|c}
S^{-1}_{(a_1,a_2,\ldots,a_d;b)}(1) & S^{-1}_{(a_1,a_2,\ldots,a_d;b)}(0) & \cdots & S^{-1}_{(a_1,a_2,\ldots,a_d;b)}(2-b) & S_b(1)\\
S^{-1}_{(a_1,a_2,\ldots,a_d;b)}(2) & S^{-1}_{(a_1,a_2,\ldots,a_d;b)}(1) & \cdots & S^{-1}_{(a_1,a_2,\ldots,a_d;b)}(3-b) & S_b(2)\\
\vdots & \vdots & \ddots & \vdots   & \vdots \\
S^{-1}_{(a_1,a_2,\ldots,a_d;b)}(b-1) & S^{-1}_{(a_1,a_2,\ldots,a_d;b)}(b-2) & \cdots & S^{-1}_{(a_1,a_2,\ldots,a_d;b)}(0) & S_b(b-1) \\
1 & 1 & \cdots & 1 & 0 
\end{array} \right)\end{align}
\end{corollary}

\proof
By the Convolution Theorem and Theorem \ref{group}, the vector 
\[
(S_{(a_1,a_2,\ldots,a_d;b)}(0),S_{(a_1,a_2,\ldots,a_d;b)}(1),\ldots,S_{(a_1,a_2,\ldots,a_d;b)}(b-1))^t
\]
is the unique solution to the system of equations
\[ \left( \begin{array}{cccc|c}
S^{-1}_{(a_1,a_2,\ldots,a_d;b)}(0) & S^{-1}_{(a_1,a_2,\ldots,a_d;b)}(-1) & \cdots & S^{-1}_{(a_1,a_2,\ldots,a_d;b)}(1-b) & S_b(0)\\
S^{-1}_{(a_1,a_2,\ldots,a_d;b)}(1) & S^{-1}_{(a_1,a_2,\ldots,a_d;b)}(0) & \cdots & S^{-1}_{(a_1,a_2,\ldots,a_d;b)}(2-b) & S_b(1)\\
\vdots & \vdots & \ddots & \vdots   & \vdots \\
S^{-1}_{(a_1,a_2,\ldots,a_d;b)}(b-1) & S^{-1}_{(a_1,a_2,\ldots,a_d;b)}(b-2) & \cdots & S^{-1}_{(a_1,a_2,\ldots,a_d;b)}(0) & S_b(b-1) \\
1 & 1 & \cdots & 1 & 0 
\end{array} \right)\]
Since $\sum_{t=0}^{b-1} S^{-1}_{(a_1,a_2,\ldots,a_d;b)}(t)=0$, the first row is equal to the sum of the second through $b$th row, times $-1$.
\proofend

\begin{example} \label{S134}
We take $d=2,a_1=1,a_2=3$ and $b=4$. In this case, Theorem \ref{systemOfEquations} states that  $S_{(1,3;4)}$ is the unique solution of the system
\[ \left( \begin{array}{cccc|c}
-1 & 2 & -1 & 0 & -\frac{1}{4} \\
0 & -1 & 2 & -1 & -\frac{1}{4} \\
-1 & 0 & -1 & 2 & -\frac{1}{4} \\
1 & 1 & 1 & 1 & 0 
\end{array} \right)\]
Solving, we see that 
\[
S_{(1,3;4)}=(\frac{5}{16},-\frac{1}{16},-\frac{3}{16},-\frac{1}{16}).
\]
\end{example}

We note that a basis for the vector space $V_b^0$ is given by
\begin{align} \label{basis}
e_i:=\delta_{\Z}(\frac{t+i}{b})-\frac{1}{b}, i=0,1,\ldots,b-2.
\end{align}

The shift operator on $V_b^0$ satisfies the relation 
\[
I+T+T^2+\ldots+T^{b-1}=0.
\]

Using the basis (\ref{basis}), we see that $I,T,T^2,\ldots,T^{b-2}$ are linearly independent. These clearly span ${\R}[T]$, and there is an isomorphism 
\begin{align} \label{iso}
{\R}[T] \simeq {\R}[X] \diagup (1+X+X^2+\ldots+X^{b-1})
\end{align}
sending $T$ to $X$.

As an application, we characterize the Fourier--Dedekind sums $S_{(a_1,a_2,\ldots,a_d;b)}$ for which 
\[
S_{(a_1,a_2,\ldots,a_d;b)}(1)=S_{(a_1,a_2,\ldots,a_d;b)}(2)=\ldots=S_{(a_1,a_2,\ldots,a_d;b)}(b-1).
\]

To give an example, we recall that the $n$th cyclotomic polynomial $\Phi_n(x)$ is defined by
\[
\Phi_n(x)= \prod_{\substack{1 \leq j \leq n,\\ (j,n)=1}} (x-\xi_n^j)
\]
and that for prime $n$, $\Phi_n(x)=x^{n-1}+x^{n-2}+\ldots+x+1$. Let $p$ be prime. Then
\[
S_{(1,2,\ldots,p-1;p)}(t)=\frac{1}{p} \sum_{j=1}^{p-1} \frac{\xi_p^{j t}}{\Phi_p(1)}=\frac{p \delta_{\Z}(\frac{t}{p})-1}{p \Phi_p(1)}=\frac{p \delta_{\Z}(\frac{t}{p})-1}{p^2}.
\]

\begin{proposition} A Fourier--Dedekind sum $S_{(a_1,a_2,\ldots,a_d;b)}(t)$ satisfies
\[
S_{(a_1,a_2,\ldots,a_d;b)}(1)=S_{(a_1,a_2,\ldots,a_d;b)}(2)=\ldots=S_{(a_1,a_2,\ldots,a_d;b)}(b-1).
\]
if and only if there exists a constant $c$ such that
\[
(1-X^{a_1})(1-X^{a_2})\cdots(1-X^{a_d}) \equiv c \pmod{1+X+X^2+\ldots+X^{b-1}}.
\]
\end{proposition}

In preparation for the proof, we state a few observations. Clearly, $T^a e_0=e_{a}$ for $a \leq b-1$. If 
\[
(c_0I+c_1T+c_2T^2+\ldots+c_{b-2}T^{b-2})e_0=0,
\]
then
\[
c_0 e_0+c_1 e_1+\ldots+c_{b-2} e_{b-2}=0,
\]
so that $c_i=0$ for each $i$.

\proof
The case $d=0$ is trivial, so we assume that $d \geq 1$. By Corollary \ref{systemOfEquations}, if 
\[
S_{(a_1,a_2,\ldots,a_d;b)}(t)=k
\]
for $t=1,2,\ldots,b-1$, then for $1 \leq i \leq b-1$,
\[
S^{-1}_{(a_1,a_2,\ldots,a_d;b)}(i) S_{(a_1,a_2,\ldots,a_d;b)}(0)+k\left( \sum_{j \neq i} S^{-1}_{(a_1,a_2,\ldots,a_d;b)}(j)\right)=-\frac{1}{b}.
\]
Since $\sum S^{-1}_{(a_1,a_2,\ldots,a_d;b)}(j)=0$ and $\sum S_{(a_1,a_2,\ldots,a_d;b)}(j)=0$, this implies that 
\[
-(b-1)k S^{-1}_{(a_1,a_2,\ldots,a_d;b)}(i)-kS^{-1}_{(a_1,a_2,\ldots,a_d;b)}(i)=-\frac{1}{b}
\]
for each $1 \leq i \leq b-1$. By definition of $S^{-1}_{(a_1,a_2,\ldots,a_d;b)}(i)$, we must have
\[
(I-T^{a_1})(I-T^{a_2})\cdots(I-T^{a_d})\delta_{\Z}(\frac{t}{b})
\]
constant for $1 \leq t \leq b-1$. Clearly,
\[
(I-T^{a_1})(I-T^{a_2})\cdots(I-T^{a_d})\delta_{\Z}(\frac{t}{b})=(I-T^{a_1})(I-T^{a_2})\cdots(I-T^{a_d})e_0,
\]
where $e_0$ is the basis vector from (\ref{basis}). In order for it to be constant on $1 \leq t \leq b-1$, we must have
\[
(I-T^{a_1})(I-T^{a_2})\cdots(I-T^{a_d})e_0=c e_0
\]
\[
\left((I-T^{a_1})(I-T^{a_2})\cdots(I-T^{a_d})-cI\right)e_0=0.
\]
By the comment above, we see that 
\[
(I-T^{a_1})(I-T^{a_2})\cdots(I-T^{a_d}) \equiv cI \pmod{I+T+T^2+\ldots+T^{b-1}}.
\]
By the isomorphism (\ref{iso}), we see that 
\[
(1-X^{a_1})(1-X^{a_2})\cdots(1-X^{a_d}) \equiv c \pmod{1+X+X^2+\ldots+X^{b-1}}.
\]
Conversely, we have
\[
S^{-1}_{(a_1,a_2,\ldots,a_d;b)}(t)=(I-T^{a_1})(I-T^{a_2})\cdots(I-T^{a_d})\delta_{\Z}(\frac{t}{b})=c\delta_{\Z}(\frac{t}{b}),
\] 
so that the converse follows from Corollary \ref{systemOfEquations}.
\proofend

\begin{lemma}
Let $p \in \N$ be prime, and let $a_1,a_2,\ldots,a_d \in \N$ satisfy $a_1 \leq a_2 \leq \ldots \leq a_d \leq p-1$ and $(a_i,p)=1$ for each $i=1,2,\ldots,d$. For $i\in \{1,2,\ldots,p-1\}$, let $e_i$ denote the multiplicity of $i$ in the multiset $\{a_1,a_2,\ldots,a_d\}$.
 A tuple $(d,a_1,a_2,\ldots,a_d,p)$ is a solution to
\[
(1-X^{a_1})(1-X^{a_2})\cdots(1-X^{a_d}) \equiv c \pmod{1+X+X^2+\ldots+X^{p-1}}
\]
for some constant $c$ if and only if
\begin{itemize} 
\item[(i).] $e:=\frac{e_1+e_2+\ldots+e_{p-1}}{p-1}=\frac{d}{p-1} \in \Z$
\item[(ii).] $e_k+e_{p-k} = 2e, \forall k \in \{1,2,\ldots,p-1\}$
\item[(iii).] $p \mid \sum_{k=1}^{\frac{p-1}{2}} k(e-e_{p-k})=-\sum_{k=1}^{\frac{p-1}{2}} k(e_k-e)$.
\end{itemize}

\end{lemma}

For the proof, we recall some facts about the $p$-th cyclotomic field for $p$ an odd prime. Let $u_k=\frac{1-\xi_p^k}{1-\xi_p}$. These units, when restricted to $2\leq k \leq \frac{p-1}{2}$ are multiplicatively independent. That is, if $c_k \in \N$ for each $k$, then
\[
\prod_{k=2}^{\frac{p-1}{2}} u_k^{c_k}= \pm \xi_p^j \implies c_k=0, \forall k.
\]
There are trivial relations for unrestricted $1 \leq k \leq p-1$. These are generated by
\begin{align} \label{trivialRelations}
u_{p-k}=- \xi^{-k} u_k.
\end{align}
We also note that if 
\[
(1-\xi)^{e_1} (1-\xi^2)^{e_2} \cdots(1-\xi^{p-1})^{e_{p^r-1}} = \pm (1-\xi)^{e_1'} (1-\xi^2)^{e_2'} \cdots(1-\xi^{p-1})^{e_{p-1}'},
\]
then by taking norms, we see that $\sum e_i = \sum e_i'$. Also,  
\[
(1-\xi)(1-\xi^2)\cdots(1-\xi^{p-1})=p.
\]

\proof

Assume first that $b=p$ for a prime $p$. Let 
\[
A=(1-\xi)(1-\xi^2)\cdots(1-\xi^{\frac{p-1}{2}}).
\]
We have
\[
(1-\xi^{\frac{p+1}{2}})(1-\xi^{\frac{p+3}{2}})\cdots(1-\xi^{p-1})=A \cdot \prod_{k=1}^{\frac{p-1}{2}} (-\xi^{-k})=(-1)^{\frac{p-1}{2}}  \xi^{-\frac{p^2-1}{8}}  A.  
\]
Using $(1-\xi)(1-\xi^2)\cdots(1-\xi^{p-1})=p$, we see that
\[
p=(-1)^{\frac{p-1}{2}}  \xi^{-\frac{p^2-1}{8}}  A^2,
\]
so that
\begin{align}
A^2= (-1)^{\frac{p-1}{2}} p \xi^{\frac{p^2-1}{8}}.
\end{align}

We evaluate both sides of

\[
(1-X^{a_1})(1-X^{a_2})\cdots(1-X^{a_d}) \equiv c \pmod{1+X+X^2+\ldots+X^{p-1}}
\]

at $X=\xi=\xi_p$:
\[
(1-\xi^{a_1})(1-\xi^{a_2})\cdots(1-\xi^{a_d})=c.
\]
Taking norms of both sides,
\[
p^d=c^{p-1} \implies c= \pm p^{\frac{d}{p-1}}.
\]
Therefore $p-1 \mid d$. Set $e=\frac{d}{p-1} \in \Z$.  We have
\[
(1-\xi^{a_1})(1-\xi^{a_2})\cdots(1-\xi^{a_d})=\pm p^e=\pm (1-\xi)^e (1-\xi^2)^e \cdots (1-\xi^{p-1})^e.
\]

We rewrite $(1-\xi^{a_1})(1-\xi^{a_2})\cdots(1-\xi^{a_d})$ as $(1-\xi)^{e_1} (1-\xi^2)^{e_2} \cdots(1-\xi^{p-1})^{e_{p-1}}$. Then 

\[
(1-\xi)^{e_1} (1-\xi^2)^{e_2} \cdots(1-\xi^{p-1})^{e_{p-1}}=\pm (1-\xi)^e (1-\xi^2)^e \cdots (1-\xi^{p-1})^e,
\]

so that $\sum e_i = (p-1)e$. Dividing both sides by $(1-\xi)^{(p-1)e}$,

\[
(\frac{1-\xi^2}{1-\xi})^{e_2} \cdots(\frac{1-\xi^{p-1}}{1-\xi})^{e_{p-1}}=\pm  (\frac{1-\xi^2}{1-\xi})^{e} \cdots(\frac{1-\xi^{p-1}}{1-\xi})^{e}
\]
implying that
\[
\mu_2^{e_2} \cdots \mu_{p-1}^{e_{p-1}} = \pm \mu_2^{e} \cdots \mu_{p-1}^{e}.
\]
Using the trivial relations,

\[
\xi^f \mu_2^{e_2+e_{p-2}} \cdots \mu_{\frac{p-1}{2}}^{e_{\frac{p-1}{2}}+e_{\frac{p+1}{2}}} = \pm \xi^{f'} \mu_2^{2e} \cdots \mu_{\frac{p-1}{2}}^{2e},
\]
for some $f,f'$.
 Consequently,
 \[
\mu_2^{e_2+e_{p-2}-2e} \cdots \mu_{\frac{p-1}{2}}^{e_{\frac{p-1}{2}}+e_{\frac{p+1}{2}}-2e} = \pm \xi^{f'-f},
\]
which implies that $e_k+e_{p-k}=2e$ for each $2 \leq k \leq \frac{p-1}{2}$. By symmetry it also follows that $e_1+e_{p-1}=2e$.

So far, we have shown that a necessary condition is that $e_k+e_{p-k}=2e$, where $e=\frac{d}{p-1} \in \N$. Assuming these conditions, we have

\[
(1-\xi)^{e_1} \cdots (1-\xi^{p-1})^{e_{p-1}}=(1-\xi)^{2e} \cdots (1-\xi^{\frac{p-1}{2}})^{2e} \prod_{k=1}^{\frac{p-1}{2}} (- \xi^{-k})^{e_{p-k}}
\]
\[
=A^{2e} \prod_{k=1}^{\frac{p-1}{2}} (- \xi^{-k})^{e_{p-k}}=\pm p^e \xi^{\frac{e(p^2-1)}{8}} \prod_{k=1}^{\frac{p-1}{2}} \xi^{-k e_{p-k}}.
\]
This expression is an integer exactly when
\[
p \mid \frac{e(p^2-1)}{8}-\sum_{k=1}^{\frac{p-1}{2}} k e_{p-k}=\sum_{k=1}^{\frac{p-1}{2}} k(e-e_{p-k}).
\]
This completes the proof in the prime case.
\proofend

\begin{corollary}
Let $p$ be prime and let $a_1,a_2,\ldots,a_d \in \N$ satisfy $a_1 \leq a_2 \leq \ldots \leq a_d \leq p-1$. For $i\in \{1,2,\ldots,p-1\}$, let $e_i$ denote the multiplicity of $i$ in the multiset $\{a_1,a_2,\ldots,a_d\}$. A Fourier--Dedekind sum $S_{(a_1,a_2,\ldots,a_d;p)}(t)$ satisfies
\[
S_{(a_1,a_2,\ldots,a_d;b)}(1)=S_{(a_1,a_2,\ldots,a_d;b)}(2)=\ldots=S_{(a_1,a_2,\ldots,a_d;b)}(p-1).
\]
if and only if 
\begin{itemize} 
\item[(i).] $e:=\frac{e_1+e_2+\ldots+e_{p-1}}{p-1}=\frac{d}{p-1} \in \Z$
\item[(ii).] $e_k+e_{p-k} = 2e, \forall k \in \{1,2,\ldots,p-1\}$
\item[(iii).] $p \mid \sum_{k=1}^{\frac{p-1}{2}} k(e-e_{p-k})=-\sum_{k=1}^{\frac{p-1}{2}} k(e_k-e)$.
\end{itemize}
\end{corollary}

Next, we show the following surprising result: the determinant of the matrix (\ref{matrix}) is equal to $(-1)^{b-1} b^d$. In particular, it is independent of the choice of $a_1,a_2,\ldots,a_d$. It is also interesting to note that if $a_i$ are not relatively prime to $b$, then the determinant is zero.

\begin{proposition} \label{determinantThm} For $d \geq 1$ and $b \geq 2$,
\begin{align} \label{det}  \begin{vmatrix}
S^{-1}_{(a_1,a_2,\ldots,a_d;b)}(1) & S^{-1}_{(a_1,a_2,\ldots,a_d;b)}(0) & \cdots & S^{-1}_{(a_1,a_2,\ldots,a_d;b)}(2-b) \\
S^{-1}_{(a_1,a_2,\ldots,a_d;b)}(2) & S^{-1}_{(a_1,a_2,\ldots,a_d;b)}(1) & \cdots & S^{-1}_{(a_1,a_2,\ldots,a_d;b)}(3-b) \\
\vdots & \vdots & \ddots & \vdots   \\
S^{-1}_{(a_1,a_2,\ldots,a_d;b)}(b-1) & S^{-1}_{(a_1,a_2,\ldots,a_d;b)}(b-2) & \cdots & S^{-1}_{(a_1,a_2,\ldots,a_d;b)}(0)  \\
1 & 1 & \cdots & 1 
\end{vmatrix} =(-1)^{b-1} b^d.
\end{align}
\end{proposition}

For what follows, we denote by $[x]$ the representative of $x+b{\Z}$ which lies in $(0,b]$.
The main idea in the proof will be to reduce the calculation of the determinant corresponding to $S^{-1}_{(a_1,a_2,\ldots,a_d;b)}$ to that of $S^{-1}_{(a_2,\ldots,a_d;b)}$.

\proof

Set
\[ M_{(a_1,a_2,\ldots,a_d;b)}:=\left(\begin{array}{cccc}
S^{-1}_{(a_1,a_2,\ldots,a_d;b)}(1) & S^{-1}_{(a_1,a_2,\ldots,a_d;b)}(0) & \cdots & S^{-1}_{(a_1,a_2,\ldots,a_d;b)}(2-b) \\
S^{-1}_{(a_1,a_2,\ldots,a_d;b)}(2) & S^{-1}_{(a_1,a_2,\ldots,a_d;b)}(1) & \cdots & S^{-1}_{(a_1,a_2,\ldots,a_d;b)}(3-b) \\
\vdots & \vdots & \ddots & \vdots   \\
S^{-1}_{(a_1,a_2,\ldots,a_d;b)}(b-1) & S^{-1}_{(a_1,a_2,\ldots,a_d;b)}(b-2) & \cdots & S^{-1}_{(a_1,a_2,\ldots,a_d;b)}(0)  \\
1 & 1 & \cdots & 1 
\end{array}\right)
\]
 Observe that
\begin{align} \label{diffDeltas}
S^{-1}_{(a_1,a_2,\ldots,a_d;b)}(t)=(I-T^{a_1})S^{-1}_{(a_2,\ldots,a_d;b)}(t).
\end{align}
Thus the matrix  $M_{(a_1,a_2,\ldots,a_d;b)}$ corresponding to $S^{-1}_{(a_1,a_2,\ldots,a_d;b)}$ is \textit{almost} obtained by replacing the $i^{th}$ row of $M_{(a_2,\ldots,a_d;b)}$ with its $i^{th}$ row minus its $[i+a_1]^{th}$ row. The issue is that the row of $1$'s prevents the matrix from being circulant, i.e., breaks the pattern. Nevertheless, we can modify the matrix corresponding to the operation of subtracting rows through use of the relation $\sum_{t} S^{-1}_{(a_2,\ldots,a_d;b)}(t)=0$.

Without loss of generality, we assume that $0 < a_i < b$ for each $i$. Define a $b \times b$ matrix $R$ as follows. Let the $i^{th}$ row $R_i$ of $R$ be
\[
R_i=\begin{cases} (e_1+e_2+\ldots+e_{b})^t+(e_i-e_{[i+a_1]})^t & \mbox{if } i=[-a_1] \\
e_b^t & \mbox{if } i=b \\
(e_i-e_{[i+a_1]})^t & \mbox{otherwise.} \\
 \end{cases}
\]

We now show that $M_{(a_1,a_2,\ldots,a_d;b)} = R M_{(a_2,\ldots,a_d;b)}$. Consider the rows of $R M_{(a_2,\ldots,a_d;b)}$. For each $i$ such that $i \neq b,[-a_1]$, the $i^{th}$ row of $R M_{(a_2,\ldots,a_d;b)}$ is equal to the difference of row $i$ of $M_{(a_2,\ldots,a_d;b)}$ and row $i+[a_1]$. By (\ref{diffDeltas}), row $i$ of $R M_{(a_2,\ldots,a_d;b)}$ agrees with row $i$ of $M_{(a_1,a_2,\ldots,a_d;b)}$.
If $i=b$, then $R_i=e_b^t$, so row $b$ of $R M_{(a_2,\ldots,a_d;b)}$ is $(e_1+e_2+\ldots+e_b)^t$, which is is also row $b$ of $M_{(a_1,a_2,\ldots,a_d;b)}$. It remains to see that rows $[-a_1]$ are equal. Row $[-a_1]$ of $R M_{(a_2,\ldots,a_d;b)}$ is equal to the sum $\Sigma$ of all rows of  $M_{(a_2,\ldots,a_d;b)}$ except the row of $1$'s, plus its $[-a_1]^{th}$ row. Since $\sum_{t \neq t'} S^{-1}_{(a_2,\ldots,a_d;b)}(t)=-S^{-1}_{(a_2,\ldots,a_d;b)}(t')$,
\[
\Sigma=-(S^{-1}_{(a_2,\ldots,a_d;b)}(0),S^{-1}_{(a_2,\ldots,a_d;b)}(b-1),\ldots,S^{-1}_{(a_2,\ldots,a_d;b)}(1)).
\]
Thus (\ref{diffDeltas}) implies that $M_{(a_1,a_2,\ldots,a_d;b)}=R M_{(a_2,\ldots,a_d;b)}$.

Next we calculate the determinant of $R$. Decompose $R$ as the sum of $R',R''$, where the rows are given by
\[
R'_i=\begin{cases} (e_1+e_2+\ldots+e_{b})^t & \mbox{if } i=[-a_1] \\
e_b^t & \mbox{if } i=b \\
(e_i-e_{[i+a_1]})^t & \mbox{otherwise.} \\
 \end{cases}
\]
\[
R''_i=\begin{cases} (e_i-e_{[i+a_1]})^t & \mbox{if } i=[-a_1] \\
e_b^t & \mbox{if } i=b \\
(e_i-e_{[i+a_1]})^t & \mbox{otherwise.} \\
 \end{cases}
\]
By linearity of the determinant on rows, $\det(R)=\det(R')+\det(R'')$. By repeatedly expanding across rows with a single nonzero entry $1$, it is easy to check that $\det(R'')=1$. As for, $R'$, we find its characteristic polynomial by solving for the eigenvalues. Let $\lambda \neq 1$ be an eigenvalue of $R'$. Suppose that $(x_1,x_2,\ldots,x_b)^t$ is an eigenvector. The $b^{th}$ row implies that $x_b=0$ so we disregard the entry $x_b$. For $i \neq [-a_1]$,
\begin{align} \label{lambda1}
x_i-x_{[i+a_1]}=\lambda x_i
\end{align}
and
\begin{align} \label{lambda2}
x_1+x_2+\ldots+x_b= \lambda x_{[-a_1]}.
\end{align}
From (\ref{lambda1}), we may write any $x_i$, $i \neq [-a_1]$ in the form $\frac{x_{[-a_1]}}{(1-\lambda)^{k_i}}$, with distinct $k_i \in \{-1,-2,\ldots,-(b-1)\}$. Substituting into (\ref{lambda2}), and dividing out by $x_{[-a_1]}$ (which cannot be zero) we obtain
\[
(1-\lambda)+\sum_{j=1}^{b-1} (1-\lambda)^{-j}=0.
\]
Multiplying through by $(1-\lambda)^{b-1}$ we see that $\det(R')=b-1$, so that $\det(R)=b$.

Finally, we find the determinant of the matrix $M_{(a;b)}$. Observe that the sum of the elements of each row except the last is equal to zero. Adding columns $2$ through $b$ to the first column shows that
\[
 \begin{vmatrix}
S^{-1}_{(a;b)}(1) & S^{-1}_{(a;b)}(0) & \cdots & S^{-1}_{(a;b)}(2-b) \\
S^{-1}_{(a;b)}(2) & S^{-1}_{(a;b)}(1) & \cdots & S^{-1}_{(a;b)}(3-b) \\
\vdots & \vdots & \ddots & \vdots   \\
S^{-1}_{(a;b)}(b-1) & S^{-1}_{(a;b)}(b-2) & \cdots & S^{-1}_{(a;b)}(0)  \\
1 & 1 & \cdots & 1 
\end{vmatrix} =
\] 
\[
 (-1)^{b-1} b\begin{vmatrix}
S^{-1}_{(a;b)}(0) & S^{-1}_{(a;b)}(-1) & \cdots & S^{-1}_{(a;b)}(2-b) \\
S^{-1}_{(a;b)}(1) & S^{-1}_{(a;b)}(0) & \cdots & S^{-1}_{(a;b)}(3-b) \\
\vdots & \vdots & \ddots & \vdots   \\
S^{-1}_{(a;b)}(b-2) & S^{-1}_{(a;b)}(b-3) & \cdots & S^{-1}_{(a;b)}(0)   
\end{vmatrix}.
\]

Thus the problem reduces to considering the top-right $(b-1) \times (b-1)$ submatrix.  Without loss of generality, we assume that $0 < a < b$. The $i^{th}$ row of the matrix is
\[
r_i=\begin{cases} (e_i-e_{i+a})^t & \mbox{if } i+a < b \\
e_i^t & \mbox{if } i+a=b \\
(e_i-e_{i+a-b})^t & \mbox{if } i+a>b, \end{cases}
\]
where $e_i$ is the $i^{th}$ standard basis vector.
\begin{example}
As an illustration, the submatrix for $a=2,b=5$ is:

\[
\left(\begin{array}{cccc}
1 & 0 & -1 &  0 \\
0 & 1 & 0 & -1 \\
0 & 0 & 1 & 0 \\
-1 & 0 & 0 & 1 
\end{array} \right).
\] 
\end{example}
 We add row $r_{[-a]}$ to row $r_{[-2a]}$, the new row $r_{[-2a]}$ to $r_{[-3a]}$ and so on until all $-1$'s have been eliminated. The result is a diagonal matrix with $1$'s across the diagonal, and consequently the submatrix has determinant $1$. 

\proofend

\begin{corollary}
Let 
\[ M_{(a_1,a_2,\ldots,a_d;b)}:=\left(\begin{array}{cccc}
S^{-1}_{(a_1,a_2,\ldots,a_d;b)}(1) & S^{-1}_{(a_1,a_2,\ldots,a_d;b)}(0) & \cdots & S^{-1}_{(a_1,a_2,\ldots,a_d;b)}(2-b) \\
S^{-1}_{(a_1,a_2,\ldots,a_d;b)}(2) & S^{-1}_{(a_1,a_2,\ldots,a_d;b)}(1) & \cdots & S^{-1}_{(a_1,a_2,\ldots,a_d;b)}(3-b) \\
\vdots & \vdots & \ddots & \vdots   \\
S^{-1}_{(a_1,a_2,\ldots,a_d;b)}(b-1) & S^{-1}_{(a_1,a_2,\ldots,a_d;b)}(b-2) & \cdots & S^{-1}_{(a_1,a_2,\ldots,a_d;b)}(0)  \\
1 & 1 & \cdots & 1 
\end{array}\right)
\]
and let $M^{(t)}_{(a_1,a_2,\ldots,a_d;b)}$ denote the matrix $M_{(a_1,a_2,\ldots,a_d;b)}$ in which column $t$ is replaced with the column vector $(-\frac{1}{b},-\frac{1}{b},\ldots,-\frac{1}{b},0)^t$. Then
\begin{align} \label{explicitFormula}
S_{(a_1,a_2,\ldots,a_d;b)}(t)=(-1)^{b-1} \frac{\det(M^{(t)}_{(a_1,a_2,\ldots,a_d;b)})}{b^d}
\end{align}
\end{corollary}

\proof
The result follows from Cramer's rule, Proposition \ref{systemOfEquations} and Proposition \ref{determinantThm}.
\proofend

The next result shows that a $(d+1)$-dimensional Fourier--Dedekind sum is a $\Z$-linear combination of $d$-dimensional Fourier--Dedekind sums.

\begin{theorem} \label{linearComb} Let $d \geq 1$. Then
\begin{align}
b S_{(a_1,a_2,\ldots,a_d)}(t)=- \sum_{k=1}^{b-1} k S_{(a_1,a_2,\ldots,a_{d-1})}(t+k a_d)
\end{align}
\end{theorem}

As before, we denote by $[x]$ the representative of $x+b{\Z}$ which lies in $(0,b]$.

\proof

By Lemma \ref{inverseLemma}, the vector $(S_{(a_1,\ldots,a_d)}(0),S_{(a_1,\ldots,a_d)}(1),\ldots,S_{(a_1,\ldots,a_d)}(b-1))^t$ is the unique solution to

\begin{align}  \left( \begin{array}{cccc|c}
S^{-1}_{(a_d;b)}(0) & S^{-1}_{(a_d;b)}(1) & \cdots & S^{-1}_{(a_d;b)}(b-1) & S_{(a_1,\ldots,a_{d-1})}(0)\\
S^{-1}_{(a_d;b)}(-1) & S^{-1}_{(a_d;b)}(0) & \cdots & S^{-1}_{(a_d;b)}(b-2) & S_{(a_1,\ldots,a_{d-1})}(1)\\
\vdots & \vdots & \ddots & \vdots   & \vdots \\
S^{-1}_{(a_d;b)}(1) & S^{-1}_{(a_d;b)}(2) & \cdots & S^{-1}_{(a_d;b)}(b-1) & S_{(a_1,\ldots,a_{d-1})}(b-1) \\
1 & 1 & \cdots & 1 & 0 
\end{array} \right).\end{align}

For $1 \leq i \leq b-1$, the $i^{th}$ row $r_i$ of the (unaugmented) matrix is $(e_i-e_{[i+a_d]})^t$. We add $-k r_{[t+k a_d]},k=1,2,\ldots,b-1$ to the $b^{th}$ row, $(1,1,\ldots,1)$. The  result is
\[
r_b-\sum_{k=1}^{b-1} k r_{[t+k a_d]}=(e_1+\ldots+e_b)^t-\sum_{k=1}^{b-1}k (e_{[t+ka_d]}-e_{[t+(k+1)a_d]})^t
=b e_t.
\]

\proofend

Taking Theorem \ref{linearComb} to its logical conclusion yields the following description of a Fourier--Dedekind sum.

\begin{corollary}
 \begin{align} \label{deltasum}
  S_{(a_1,a_2,\ldots,a_d;b)}(t)=
\frac{(-1)^d}{b^d}[\sum_{\substack{1 \leq k_1, k_2,\ldots, k_d \leq b-1,\\ a_1 k_1+a_2 k_2 +\ldots+a_d k_d \equiv -t \pmod{b}}} k_1 k_2 \cdots k_d-\frac{1}{b}\binom{b}{2}^d]  
\end{align}
\end{corollary}

\proof
Inductively, we see that 
\[
b^d S_{(a_1,a_2,\ldots,a_d;b)}(t)=(-1)^d \sum_{1 \leq k_1, k_2,\ldots, k_d \leq b-1} k_1 k_2 \cdots k_d S_b(t+k_1 a_1+k_2 a_2 +\ldots+k_d a_d).
\]
By Lemma \ref{geometricSeries}, $S_b(t)=\delta_{\Z}(\frac{t}{b})-\frac{1}{b}$. Simplifying yields
\[
 S_{(a_1,a_2,\ldots,a_d;b)}(t)=\frac{(-1)^d}{b^d}[\sum_{1 \leq k_1, k_2,\ldots, k_d \leq b-1} k_1 k_2 \cdots k_d \delta_{\Z}(\frac{t+k_1 a_1+k_2 a_2+\ldots+k_d a_d}{b})-\frac{1}{b}\binom{b}{2}^d].
 \]
Since
\[
\sum_{1 \leq k_1, k_2,\ldots, k_d \leq b-1} k_1 k_2 \cdots k_d \delta_{\Z}(\frac{t+k_1 a_1+k_2 a_2+\ldots+k_d a_d}{b})=\sum_{\substack{1 \leq k_1, k_2,\ldots, k_d \leq b-1,\\ a_1 k_1+a_2 k_2 +\ldots+a_d k_d \equiv -t \pmod{b}}} k_1 k_2 \cdots k_d,
\]
the result follows.
\proofend

\begin{remark}
Formula (\ref{deltasum}) allows us to generalize Fourier--Dedekind sums to arguments $a_i$ which are not necessarily relatively prime to $b$.
\end{remark}

\begin{definition}  For any $a_1,a_2,\ldots,a_d,b \in \N$, we define the Reduced Fourier--Dedekind sum $\tilde{S}_{(a_1,a_2,\ldots,a_d;b)}(t)$ by
\begin{align} \label{nontrivialPart}
\tilde{S}_{(a_1,a_2,\ldots,a_d;b)}(t):=\sum_{\substack{1 \leq k_1, k_2,\ldots, k_d \leq b-1,\\ a_1 k_1+a_2 k_2 +\ldots+a_d k_d \equiv -t \pmod{b}}} k_1 k_2 \cdots k_d.
\end{align}
\end{definition}
Conceptually, the Reduced Fourier--Dedekind sum (\ref{nontrivialPart}) is the nontrivial part of a Fourier--Dedekind sum:

\[
 S_{(a_1,a_2,\ldots,a_d;b)}(t)=
\frac{(-1)^d}{b^d}[\tilde{S}_{(a_1,a_2,\ldots,a_d;b)}(t)-\frac{1}{b}\binom{b}{2}^d].
\]

 Our aim now will be to better understand these functions.

The next result describes the Reduced Fourier--Dedekind sum in terms of a relatively simple generating function. 
In words, it says that $\tilde{S}_{(a_1,a_2,\ldots,a_d;b)}(t)$ is equal to the sum of the coefficients of the monomials of $\left(z^{a_1}+2z^{2a_1}+\ldots+(b-1)z^{(b-1)a_1}\right)\cdots \left(z^{a_d}+2z^{2a_d}+\ldots+(b-1)z^{(b-1)a_d}\right )$ which have exponent congruent to $-t$ modulo $b$.

\begin{theorem} \label{genFuncThm}
For any $a_1,a_2,\ldots,a_d,b \in \N$, \\\\
$
 \tilde{S}_{(a_1,a_2,\ldots,a_d;b)}(t)=
$
\begin{align} \label{genFunc}
 \sum_{j= -\infty}^{\infty} [z^{-t+bj}] \left(z^{a_1}+2z^{2a_1}+\ldots+(b-1)z^{(b-1)a_1}\right)\cdots \left(z^{a_d}+2z^{2a_d}+\ldots+(b-1)z^{(b-1)a_d}\right ).
 \end{align}
\end{theorem}

\proof
We have
\[
[z^{-t+bj}] \left(z^{a_1}+2z^{2a_1}+\ldots+(b-1)z^{(b-1)a_1}\right)\cdots \left(z^{a_d}+2z^{2a_d}+\ldots+(b-1)z^{(b-1)a_d}\right ).
\]
\[
=\sum_{\substack{1 \leq k_1, k_2,\ldots, k_d \leq b-1,\\ a_1 k_1+a_2 k_2 +\ldots+a_d k_d = -t+bj}} k_1 k_2 \cdots k_d.
\]
Summing over $j$ yields the result.
\proofend

 One consequence of the interpretation of Theorem \ref{genFuncThm} is that, just like Fourier--Dedekind sums, the Reduced Fourier--Dedekind sums can be built up from lower dimensional Fourier--Dedekind sums.
 
 \begin{corollary} For any $a_1,a_2,\ldots,a_d,b \in \N$,
 \[
\tilde{S}_{(a_1,a_2,\ldots,a_d)}(t)= \sum_{k=1}^{b-1} k \tilde{S}_{(a_1,a_2,\ldots,a_{d-1})}(t+k a_d)
 \]
 \end{corollary}
 
 \proof
  \[
  \sum_{j= -\infty}^{\infty} [z^{-t+bj}] \left(z^{a_1}+2z^{2a_1}+\ldots+(b-1)z^{(b-1)a_1}\right)\cdots \left(z^{a_d}+2z^{2a_d}+\ldots+(b-1)z^{(b-1)a_d}\right )
  \]
  \[
  =\sum_{m=1}^{b-1} m \sum_{j= -\infty}^{\infty} [z^{-t-ma_d+bj}] \left(z^{a_1}+2z^{2a_1}+\ldots+(b-1)z^{(b-1)a_1}\right)\cdots \left(z^{a_{d-1}}+2z^{2a_{d-1}}+\ldots+(b-1)z^{(b-1)a_{d-1}}\right )
  \]
  \[
  =\sum_{m=1}^{b-1} m \tilde{S}_{(a_1,a_2,\ldots,a_d;b)}(t+ma_d).
  \]
 \proofend

 \begin{corollary} \label{pairs} Let $[x]$ be the representative of $x+b{\Z}$ which lies in $[0,b)$. If $a_1,a_2$ are relatively prime to $b$, then
 \[
 \tilde{S}_{(a_1,a_2;b)}(t)= \sum_{k=1}^{b-1} k [a_2^{-1}(-t-ka_1)]=\sum_{k=1}^{b-1} k [a_1^{-1}(-t-ka_2)].
 \]
 As a consequence, 
  \[
 \tilde{S}_{(a_1,a_2;b)}(t) \equiv  -a_2^{-1} t \frac{(b-1)b}{2}-a_2^{-1} a_1 \frac{(b-1)b(2b-1)}{6}
 \]
\[ 
  \equiv -a_1^{-1} t \frac{(b-1)b}{2}-a_1^{-1} a_2 \frac{(b-1)b(2b-1)}{6} \pmod{b}.
 \]
 \end{corollary}
 
 \proof
 To each exponent $ka_1$ of $z^{a_1}+2z^{2a_1}+\ldots+(b-1)z^{(b-1)a_1}$, corresponds another exponent $\ell a_2$ in $z^{a_2}+2z^{2a_2}+\ldots+(b-1)z^{(b-1)a_2}$ such that
 \[
 k a_1+\ell a_2 \equiv - t \pmod{b}.
 \]
 Solving for $\ell$,
 \[
 \ell \equiv a_2^{-1}(-t-ka_1)  \pmod{b}.
 \]
 \proofend

We may recover the expression for $\tilde{S}_{(a_1,a_2,\ldots,a_d;b)}(t)$ in terms of roots of unity. Indeed, some  manipulation of series shows that 
\begin{align} \label{seriesCompact} 
 f_{(a;b)}(z):=\sum_{k=1}^{b-1} k z^{k a}=\frac{(b-1)z^{a(b+1)}-b z^{a b}+z^a}{(1-z^a)^2}.
 \end{align}
As a consequence,
\[
\tilde{S}_{(a_1,a_2,\ldots,a_d;b)}(t)= \sum_{j=0}^{b-1}\lim_{z \rightarrow \xi_b^j} z^t f_{(a_1;b)}(z) f_{(a_2;b)}(z) \cdots f_{(a_d;b)}(z).
\]
Since
\[
\lim_{z \rightarrow \xi_b^j} f_{(a;b)}(z)=\begin{cases} \binom{b}{2} &\mbox{if } j \equiv 0 \pmod{b} \\
-\frac{b}{1-\xi_b^{ja}} & \mbox{if } j \not \equiv 0 \pmod{b} \end{cases}.
\]
we are able to recover the definition in terms of roots of  unity.

\begin{figure}[H]
\centering
\includegraphics[width=5in]{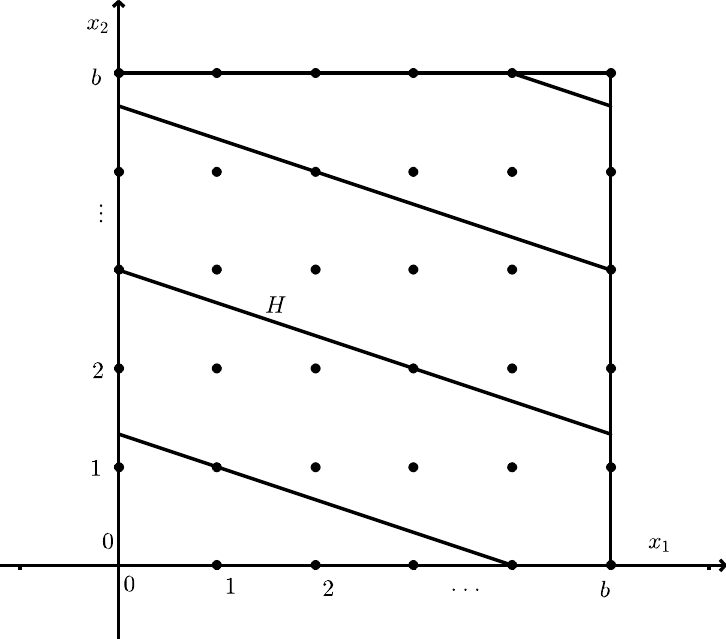}
\caption{\label{squareLattice} The torus $T$ and the hyperplane $H$ in 2 dimensions. Each lattice point on $H$ is assigned a weight $x_1 x_2$ and the sum of the weights of these lattice points is equal to $\tilde{S}_{(a_1,a_2,\ldots,a_d;b)}(t)$. }
\end{figure}

It is interesting to interpret the Reduced Fourier--Dedekind sum geometrically. 
We consider the following torus $T$. Its fundamental domain is $F=\{(x_1,x_2,\ldots,x_n) \in \R^n: 0 \leq x_1,x_2,\ldots,x_n < b\}$. We assign each point $(x_1,x_2,\ldots,x_n)$ in $F$ the weight $x_1  x_2 \cdots x_n$, which one could interpret as a suitable volume, and extend periodically to $T$. Let $H$ be the hyperplane $a_1x_1+a_2x_2+\ldots+a_nx_n \equiv -t \pmod{b}$. The value  $\tilde{S}_{(a_1,a_2,\ldots,a_d;b)}(t)$ is equal to the weighted sum over lattice points in $T \cap H$ (see Figure \ref{squareLattice} for the 2-dimensional case).

 In each row $y=i \in \Z \cap [0,b-1]$ and each column $x=j \in \Z \cap [0,b-1]$ of the fundamental domain $F$, exactly one lattice point lies on $L$. The inequalities on $\tilde{S}_{(a_1,a_2;b)}(0)$ become evident from this interpretation, and it is easy to see that the largest value of $\tilde{S}_{(a_1,a_2;b)}(0)$ occurs when the lattice points occupy the diagonal $x=y$, and the smallest when they occupy the diagonal $x=-y$.
 
 We also note that by summing the lattice points in different orders, we obtain different expressions for $\tilde{S}_{(a_1,a_2;b)}(t)$. For example, summing by rows or columns we obtain Corollary \ref{pairs} and its symmetric counterpart. An alternative method for summation is to unfold the torus.

\section{An extension of Rademacher reciprocity}

We now prove a formula relating $S_{(a_1,a_2,\ldots,a_d;b)}(-t)$ and $S_{(a_1,a_2,\ldots,a_d;b)}(t+a_1+a_2+\ldots+a_d)$ to $S_{(a_1,a_2,\ldots,a_d;b)}(t)$.

\begin{lemma} \label{PIE}
Let $a_1,a_2,\ldots,a_d,b \in \N$ with $a_i$ relatively prime to $b$ for $i=1,2,\ldots,d$. For every $t \in \Z$,
\begin{itemize}
\item[(i).] $S_{(a_1,a_2,\ldots,a_d;b)}(-t)=\sum_{k=0}^{d} (-1)^k \sum_{1 \leq i_1 < \ldots < i_k \leq d} S_{(a_{i_1},\ldots,a_{i_k};b)}(t)$.
\item[(ii).] $S_{(a_1,a_2,\ldots,a_d;b)}(t+a_1+a_2+\ldots+a_d)=\sum_{k=0}^{d} (-1)^{d-k} \sum_{1 \leq i_1 < \ldots < i_k \leq d} S_{(a_{i_1},\ldots,a_{i_k};b)}(t)$.
\end{itemize}
\end{lemma}

\proof
\begin{itemize} 
\item[(i).] Recall that $S_{(a;b)}(t)+S_{(a;b)}(-t)=S_b(t)=S_b(-t)$. By Theorem \ref{functionalEqConv},
\[
S_{(a_1,a_2,\ldots,a_d;b)}=S_{(a_1;b)} * S_{(a_2;b)} * \ldots * S_{(a_d;b)}.
\]
Note that for $b$-periodic functions $f_1, f_2, \ldots, f_m$,
\[
(f_1(x) * f_2(x) * \ldots * f_m(x))(-t)=\sum_{\substack{x_1 \bmod{b},x_2 \bmod{b},\ldots,x_m \bmod{b},\\ x_1+x_2+\ldots+x_m \equiv -t \bmod{b}}} f(x_1) f(x_2) \cdots f(x_m)
\]
\[
=\sum_{\substack{x_1 \bmod{b},x_2 \bmod{b},\ldots,x_m \bmod{b},\\ x_1+x_2+\ldots+x_m \equiv t \bmod{b}}} f(-x_1) f(-x_2) \cdots f(-x_m)
\]
\[
=(f_1(-x) * f_2(-x) * \ldots * f_m(-x))(t).
\]
Since $S_b(x)=\delta_{\Z}(\frac{x}{b})-\frac{1}{b}$, 
\[
S_b * S_b = S_b
\]
and
\[
S_b * S_{(a_1,a_2,\ldots,a_k)}=S_{(a_1,a_2,\ldots,a_k)}.
\]
Therefore
\[
S_{(a_1,a_2,\ldots,a_d;b)}(-t)=[(S_b(x)-S_{(a_1;b)}(x)) * (S_b(x)-S_{(a_2;b)}(x)) * \ldots * (S_b(x)-S_{(a_d;b)}(x))](t).
\]
\item[(ii).] We prove more generally that for $1 \leq k \leq d$,
\[
S_{(a_1,a_2,\ldots,a_d;b)}(t+a_k+a_{k+1}+\ldots+a_d)=\sum_{j=0}^{d-k+1} (-1)^{j} \sum_{k \leq i_1 <\ldots < i_j \leq d} S_{(a_1,\ldots,\hat{a}_{i_1},\hat{a}_{i_2},\ldots,\hat{a}_{i_j},\ldots,a_d;b)}(t).
\]
Proceed by induction on $1 \leq k \leq d$. By Theorem \ref{functionalEqConv},
\[
S_{(a_1,a_2,\ldots,a_d;b)}(t+a_1+a_2+\ldots+a_d)=
\]
\[
=S_{(a_1,a_2,\ldots,a_d;b)}(t+a_1+a_2+\ldots+a_{d-1})-S_{(a_1,a_2,\ldots,a_{d-1};b)}(t+a_1+a_2+\ldots+a_{d-1})
\]

Consequently, the base case $k=d=1$ holds:
\[
S_{(a_1;b)}(t+a_1)=S_{(a_1;b)}(t)-S_{b}(t),
\] 
while
\[
\sum_{j=0}^{d-k+1} (-1)^{j} \sum_{k \leq i_1 <\ldots < i_j \leq d} S_{(a_1,\ldots,\hat{a}_{i_1},\hat{a}_{i_2},\ldots,\hat{a}_{i_j},\ldots,a_d;b)}(t)=S_{(a_1;b)}(t)-S_{b}(t).
\]
Assume the inductive hypothesis.
\[
S_{(a_1,a_2,\ldots,a_d;b)}(t+a_1+a_2+\ldots+a_d)=
\]
\[
=S_{(a_1,a_2,\ldots,a_d;b)}(t+a_2+\ldots+a_d)-S_{(a_2,\ldots,a_d;b)}(t+a_2+\ldots+a_d)
\]
\[
=\sum_{j=0}^{d-1} (-1)^{j} \sum_{2 \leq i_1 <\ldots < i_j \leq d} S_{(a_1,\ldots,\hat{a}_{i_1},\hat{a}_{i_2},\ldots,\hat{a}_{i_j},\ldots,a_d;b)}(t)
\]
\[
-\sum_{j=0}^{d-1} (-1)^{j} \sum_{2 \leq i_1 <\ldots < i_j \leq d} S_{(a_2,\ldots,\hat{a}_{i_1},\hat{a}_{i_2},\ldots,\hat{a}_{i_j},\ldots,a_d;b)}(t)
\]
\[
=\sum_{j=0}^{d} (-1)^{j} \sum_{1 \leq i_1 <\ldots < i_j \leq d} S_{(a_1,\ldots,\hat{a}_{i_1},\hat{a}_{i_2},\ldots,\hat{a}_{i_j},\ldots,a_d;b)}(t).
\]
\end{itemize}
\proofend

Recall that
\[
R_{\{a_1,a_2,\ldots,a_d\}}(t):=\sum_{m=1}^{d} S_{(a_1,\ldots,\hat{a}_m,\ldots,a_d;a_m)}(t).
\]

\begin{corollary} \label{RPIE}
Let $a_1, a_2,\ldots,a_d \in \N$ be pairwise relatively prime. For all $t \in \Z$, 
\begin{itemize} 
\item[(i).] $R_{\{a_1,a_2,\ldots,a_d\}}(-t)=\sum_{k=0}^{d-1} (-1)^{k} \sum_{1 \leq i_1 < \ldots < i_{k+1} \leq d} R_{\{a_{i_1},\ldots,a_{i_{k+1}}\}}(t)$.
\item[(ii).] $R_{\{a_1,a_2,\ldots,a_d\}}(t+a_1+a_2+\ldots+a_d)=\sum_{k=0}^{d-1} (-1)^{d+k} \sum_{1 \leq i_1 < \ldots < i_{k+1} \leq d} R_{\{a_{i_1},\ldots,a_{i_{k+1}}\}}(t)$.
\end{itemize}
\end{corollary}

\proof
\begin{itemize} 
\item[(i).]  $\sum_{m=1}^{d} S_{(a_1,\ldots,\hat{a}_m,\ldots,a_d;a_m)}(-t)=\sum_{m=1}^{d} \sum_{k=0}^{d-1} (-1)^k\sum_{\substack{1 \leq i_1 < \ldots < i_k \leq d,\\ i_1,i_2,\ldots,i_k \neq m}} S_{(a_{i_1},\ldots,a_{i_k};a_m)}(t)$
\[
=\sum_{k=0}^{d-1} (-1)^k \sum_{m=1}^{d} \sum_{\substack{1 \leq i_1 < \ldots < i_k \leq d,\\ i_1,i_2,\ldots,i_k \neq m}}  S_{(a_{i_1},\ldots,a_{i_k};a_m)}(t)
\]
\[
=\sum_{k=0}^{d-1} (-1)^{k} \sum_{1 \leq i_1 < \ldots < i_{k+1} \leq d} R_{\{a_{i_1},\ldots,a_{i_{k+1}}\}}(t).
\]
\item [(ii)] The proof is analogous to the one in (i).
\end{itemize}
\proofend

\begin{lemma} \label{polynomialPIE}
Let $a_1,a_2,\ldots,a_d \in \N$ be pairwise relatively prime. Then
\[
\text{\textnormal{poly}}_{\{-a_1,-a_2,\ldots,-a_d\}}(t)=\sum_{k=1}^{d} (-1)^k \sum_{1 \leq i_1 < \ldots < i_k \leq d} \text{\textnormal{poly}}_{\{a_{i_1},\ldots,a_{i_k}\}}(t).
\]
\end{lemma}

\proof
If $h(z)=\frac{1}{(1-z^{a_1})(1-z^{a_2})\cdots(1-z^{a_d})z^n}$ and the partial fraction expansion of $h(z)$ is
\[
h(z)=\frac{A_1}{z}+\frac{A_2}{z}+\ldots+\frac{A_n}{z^n}+\frac{B_1}{z-1}+\frac{B_2}{(z-1)^2}+\ldots+\frac{B_d}{(z-1)^d}
\]
\[
+\sum_{k=1}^{a_1-1}\frac{C_{1k}}{z-\xi_{a_1}^{k}}+\sum_{k=1}^{a_2-1}\frac{C_{2k}}{z-\xi_{a_2}^k}+\ldots+\sum_{k=1}^{a_d-1}\frac{C_{dk}}{z-\xi_{a_d}^k},
\]
then
\[
\text{poly}_{\{a_1,a_2,\ldots,a_d\}}(n):=-B_1+B_2-\ldots+(-1)^d B_d.
\]

We have
\[
\frac{1}{z^n}[\sum_{k=1}^d (-1)^{k+1} \sum_{1 \leq i_1 < \ldots < i_k \leq d} \frac{1}{(1-z^{a_{i_1}})\cdots(1-z^{a_{i_k}})}]=
\]
\[
=\frac{\sum_{j=0}^{d-1} (-1)^{d-j+1} \sum_{1 \leq i_1 < \ldots < i_j \leq d} (1-z^{a_{i_1}})(1-z^{a_{i_2}})\cdots(1-z^{a_{i_j}}) }{(1-z^{a_1})(1-z^{a_2})\cdots (1-z^{a_d})z^n}.
\]
On the other hand,
\[
(-1)^{d+1} z^{a_1} z^{a_2} \cdots z^{a_d}=(-1)^{d+1} (1-(1-z^{a_1})) (1-(1-z^{a_2})) \cdots (1-(1-z^{a_d}))=
\]
\[
=- (1-z^{a_1})(1-z^{a_2})\cdots(1-z^{a_d})+\sum_{j=0}^{d-1} (-1)^{d+j+1} \sum_{1 \leq i_1 < \ldots < i_j \leq d} (1-z^{a_{i_1}})(1-z^{a_{i_2}})\cdots(1-z^{a_{i_j}}) .
\]
It follows that
\[
\frac{1}{z^n}[\sum_{k=1}^d (-1)^{k+1} \sum_{1 \leq i_1 < \ldots < i_k \leq d} \frac{1}{(1-z^{a_{i_1}})\cdots(1-z^{a_{i_k}})}]=
\]
\[
=\frac{(-1)^{d+1} z^{a_1} z^{a_2} \cdots z^{a_d}+(1-z^{a_1})(1-z^{a_2})\cdots (1-z^{a_d})}{z^n (1-z^{a_1})(1-z^{a_2})\cdots (1-z^{a_d})}
\]
\[
=\frac{1}{z^n} +\frac{(-1)^{d+1} z^{a_1 + a_2 +\ldots +a_d} }{z^n (1-z^{a_1})(1-z^{a_2})\cdots (1-z^{a_d})}.
\]
For $n$ satisfying $n+a_1+a_2+\ldots+a_d > a_1 a_2 \cdots a_d$, the last expression is a proper rational function, and can consequently be expanded into partial fractions. We recognize that the corresponding polynomial is $(-1)^{d+1} \text{poly}_{\{a_1,a_2,\ldots,a_d\}}(n-a_1- a_2 \ldots -a_d)$. Thus we have shown that for all sufficiently large $n$, the two polynomials, 
\[
(-1)^d \text{poly}_{\{a_1,a_2,\ldots,a_d\}}(n-a_1-a_2-\ldots-a_d)
\]
and 
\[
\sum_{k=1}^{d} (-1)^k \sum_{1 \leq i_1 < \ldots < i_k \leq d} \text{poly}_{\{a_{i_1},\ldots,a_{i_k}\}}(n)
\]
 take on the same values. Consequently, they must be equal for all $n$.
\proofend

We note that
\[
\text{poly}_{\{-a_1,-a_2,\ldots,-a_d\}}(t)=(-1)^d \text{poly}_{\{a_1,a_2,\ldots,a_d\}}(t-a_1-a_2-\ldots-a_d).
\]
Indeed, 
\[
\frac{(-1)^{d} z^{a_1 + a_2 +\ldots +a_d} }{z^n (1-z^{a_1})(1-z^{a_2})\cdots (1-z^{a_d})}=\frac{1 }{z^n (1-z^{-a_1})(1-z^{-a_2})\cdots (1-z^{-a_d})}.
\]
\begin{example}
 The first few instances of $\text{\textnormal{poly}}$ are \cite[Example 8.3]{BR1}
\[
\text{\textnormal{poly}}_{\{a_1\}}(t)=\frac{1}{a_1},
\]
\[
\text{\textnormal{poly}}_{\{a_1,a_2\}}(t)=\frac{t}{a_1 a_2}+\frac{1}{2}(\frac{1}{a_1}+\frac{1}{a_2})
\]
The case $d=1$ in Lemma \ref{polynomialPIE} states that $\text{\textnormal{poly}}_{\{-a_1\}}(t)=-\text{\textnormal{poly}}_{\{a_1\}}(t)$, which is clear by inspection. \\
The case $d=2$ in Lemma \ref{polynomialPIE} states that
\[
\text{\textnormal{poly}}_{\{-a_1,-a_2\}}(t)= \frac{t}{a_1 a_2}-\frac{1}{2}(\frac{1}{a_1}+\frac{1}{a_2})
\]
is equal to
\[
-\text{\textnormal{poly}}_{\{a_1\}}(t)-\text{\textnormal{poly}}_{\{a_2\}}(t)+\text{\textnormal{poly}}_{\{a_1,a_2\}}(t)=-\frac{1}{a_1}-\frac{1}{a_2}+\frac{t}{a_1 a_2}+\frac{1}{2}(\frac{1}{a_1}+\frac{1}{a_2}).
\]
\end{example}

Recall from the introduction that an explicit formula for $\text{poly}$ is given by
\[
\text{poly}_{\{a_1,a_2,\ldots,a_d\}}(t)=\frac{1}{a_1 \cdots a_d} \sum_{m=0}^{d-1} \frac{(-1)^m}{(d-1-m)!} \sum_{k_1+\ldots+k_d=m} a_1^{k_1} \cdots a_d^{k_d} \frac{B_{k_1} \cdots B_{k_d}}{k_1! \cdots k_d !}t^{d-1-m},
\]
where $B_j$ denotes the $j$th Bernoulli Number. Evidently, $\text{poly}$ is homogenous of degree $-1$ when viewed as a function of $a_1,a_2,\ldots,a_d,n$. Consequently, $\text{poly}_{\{-a_1,-a_2,\ldots,-a_d\}}(t)= \text{poly}_{\{a_1,a_2,\ldots,a_d\}}(-t)$.

The following Theorem extends Rademacher reciprocity.

\begin{theorem} \label{RademacherExtended}Let $a_1,a_2,\ldots,a_d \in \N$ be pairwise relatively prime. Let $n \in \Z$. If one of $(i)-(iii)$ holds, where
\begin{itemize}
\item[(i).] $1-\min\{a_1,a_2,\ldots,a_d\} \leq n \leq -1$
\item[(ii).] $1 \leq n \leq a_1+a_2+\ldots+a_d-1$
\item[(iii).] $a_1+a_2+\ldots+a_d+1 \leq n \leq a_1+a_2+\ldots+a_d+\min\{a_1,a_2,\ldots,a_d\}-1$
\end{itemize}
then
\[
S_{(a_2,\ldots,a_{d};a_1)}(n)+S_{(a_1,a_3,a_4,\ldots,a_{d};a_2)}(n)+\ldots+S_{(a_1,a_2,\ldots,a_{d-1};a_d)}(n)=-\text{\textnormal{poly}}_{\{a_1,a_2,\ldots,a_d\}}(-n).
\]
\end{theorem}

In the notation introduced in this paper, the Theorem states that under the assumptions of $(i),(ii)$ or $(iii)$,
\[
R_{\{a_1,a_2,\ldots,a_d\}}(n)=-\text{poly}_{\{a_1,a_2,\ldots,a_d\}}(-n).
\]

\proof
\begin{itemize} \item[(i).] By Corollary \ref{RPIE}, 
\[
R_{\{a_1,a_2,\ldots,a_d\}}(t)=\sum_{k=0}^{d-1} (-1)^{k} \sum_{1 \leq i_1 < \ldots < i_{k+1} \leq d} R_{\{a_{i_1},\ldots,a_{i_{k+1}}\}}(-t).
\]
For $1-\min\{a_1,a_2,\ldots,a_d\} \leq t \leq -1$ an integer, we may apply Rademacher reciprocity to each $R_{\{a_{i_1},\ldots,a_{i_k}\}}(-t)$:
\[
\sum_{k=0}^{d-1} (-1)^{k} \sum_{1 \leq i_1 < \ldots < i_{k+1} \leq d} R_{\{a_{i_1},\ldots,a_{i_{k+1}}\}}(-t)=\sum_{k=0}^{d-1} (-1)^{k+1} \sum_{1 \leq i_1 < \ldots < i_{k+1} \leq d} \text{poly}_{\{a_{i_1},\ldots,a_{i_{k+1}}\}}(t).
\]

By Lemma \ref{polynomialPIE},
\[
\sum_{k=0}^{d-1} (-1)^{k+1} \sum_{1 \leq i_1 < \ldots < i_{k+1} \leq d} \text{poly}_{\{a_{i_1},\ldots,a_{i_{k+1}}\}}(t)=\text{poly}_{\{-a_1,-a_2,\ldots,-a_d\}}(t)=-\text{poly}_{\{a_1,a_2,\ldots,a_d\}}(-t).
\]

\item[(ii).] This is a restatement of Rademacher reciprocity.

\item[(iii).] Let $a_1+a_2+\ldots+a_d+1 \leq n \leq a_1+a_2+\ldots+a_d+\min\{a_1,a_2,\ldots,a_d\}-1$ and set $t=n-a_1-a_2-\ldots-a_d$. 
\[
R_{\{a_1,a_2,\ldots,a_d\}}(n)=R_{\{a_1,a_2,\ldots,a_d\}}(t+a_1+a_2+\ldots+a_d).
\]
By Corollary \ref{RPIE},
\[
R_{\{a_1,a_2,\ldots,a_d\}}(t+a_1+a_2+\ldots+a_d)=\sum_{k=0}^{d-1} (-1)^{k} \sum_{1 \leq i_1 < \ldots < i_{k+1} \leq d} R_{\{a_{i_1},\ldots,a_{i_{k+1}}\}}(t).
\]
Applying Reciprocity yields the result.
\end{itemize}

\proofend

\section{Average behavior of Fourier--Dedekind sums}

 In this section, we study the average behavior of a Fourier--Dedekind sum as the $a_i$'s vary.  

\begin{definition} \label{avgdef}
The average over the $i$th variable of $S_{(a_1,a_2,\ldots,a_{d-1},a_d;b)}$ at $t$, denoted by $S_{(a_1,a_2,\ldots,\bar{a}_i,\ldots,a_d;b)}(t)$, is defined to be
\[
S_{(a_1,a_2,\ldots,\bar{a}_i,\ldots,a_d;b)}(t):=\frac{1}{\phi(b)} \sum_{\substack{1 \leq m \leq b-1\\ (m,b)=1}}  S_{(a_1,a_2,\ldots,a_{i-1},m,a_{i+1},\ldots,a_d;b)}(t).
\] 
The average over all variables of $S_{(a_1,a_2,\ldots,a_{d-1},a_d;b)}$ at $t$, denoted by $S_{(\bar{a}_1,\bar{a}_2,\ldots,\bar{a}_d;b)}(t)$, is defined to be
\[
S_{(\bar{a}_1,\bar{a}_2,\ldots,\bar{a}_d;b)}(t):=(\frac{1}{\phi(b)})^d \sum_{\substack{1 \leq m_1,m_2,\ldots,m_d \leq b-1\\ (m_i,b)=1}}  S_{(m_1,m_2,\ldots,m_d;b)}(t)
\] 
\end{definition}

\begin{theorem} \label{avg} Let $b \geq 3$ and let $(a_i,b)=1$ for each $i$. For every $t \in \Z$,
\[
S_{(a_1,a_2,\ldots,a_{d-1},\bar{a}_d;b)}(t)=\frac{1}{2} S_{(a_1,a_2,\ldots,a_{d-1})}(t).
\] 
and
\[
S_{(\bar{a}_1,\bar{a}_2,\ldots,\bar{a}_d;b)}(t)=\frac{1}{2^d} S_{b}(t)=\frac{\delta_{\Z}(\frac{t}{b})-\frac{1}{b}}{2^d} .
\] 
\end{theorem}

Before proving the result, we deduce a helpful Lemma.

\begin{lemma} \label{-a1}For every $t \in \Z$,
\[
S_{(-a_1,a_2,\ldots,a_d;b)}(t)=-S_{(a_1,a_2,\ldots,a_d;b)}(t)+S_{(a_2,\ldots,a_d;b)}(t).
\]
\end{lemma}

\proof

By \cite[Pg. 144]{BR1},
\[
S_{(a;b)}(t)=\Big( \Big(\frac{a^{-1} t}{b} \Big) \Big)+\frac{1}{2} \delta_{\Z}(\frac{t}{b})-\frac{1}{2b},
\]

where $a^{-1}$ is the modular inverse of $a$ modulo $b$.
Consequently,
\[
S_{(-a;b)}(t)=-S_{(a;b)}(t)+S_b(t).
\]

By Theorem \ref{functionalEqConv},
\[
S_{(-a_1,a_2,\ldots,a_d;b)}= S_{(-a_1;b)} * S_{(a_2;b)} * \ldots * S_{(a_d;b)}
\]
\[
=(-S_{(a_1;b)}(t)+S_b(t)) * S_{(a_2;b)} * \ldots * S_{(a_d;b)}=-S_{(a_1,a_2,\ldots,a_d;b)}+S_{(a_2,\ldots,a_d;b)}.
\]
\proofend

\noindent{\bf Proof of Theorem \ref{avg}.}
We add pairs of the form $S_{(a_1,a_2,\ldots,a_{d-1},m)}(t)$ and $S_{(a_1,a_2,\ldots,a_{d-1},-m)}(t)$ and apply Lemma \ref{-a1}. Since $b \geq 3$, each such pair represents two distinct summands.
\proofend

\section{Bounds, maxima and minima of 2-dimensional Fourier--Dedekind sums}
In this section, our aim is three-fold: we aim to obtain bounds of $2$-dimensional Fourier--Dedekind sums for a fixed $t$, to better understand the location of maxima and minima of Fourier--Dedekind sums as $t$ varies, and to find bounds on reciprocal sums of $2$-dimensional Fourier--Dedekind sums.

\begin{theorem} \label{DedekindSumBounds} For all $a_1,a_2$ coprime to $b$,

\[
-\frac{(b-1)(b-5)}{12b} \leq S_{(a_1,a_2;b)}(0) \leq \frac{(b-1)(b+1)}{12 b}.
\]
The upper bound holds if and only if $a_1+a_2 \equiv 0 \pmod{b}$. The lower bound holds if and only if $a_1 \equiv a_2 \pmod{b}$.

For all $a_1,a_2$ coprime to $b$ and $1 \leq t \leq b-1$,
\[
-\frac{(b-1)(b+1)}{12b} \leq S_{(a_1,a_2;b)}(t) \leq \frac{(b-1)(b-5)}{12 b}.
\]
The upper bound holds if and only if $a_1 \equiv -a_2 \equiv t \pmod{b}$. The lower bound holds if and only if $a_1 \equiv a_2 \equiv t \pmod{b}$.

\end{theorem}

\proof
The first part is equivalent to showing that
\[
\frac{b(b-1)(b+1)}{6} \leq \tilde{S}_{(a_1,a_2;b)}(0) \leq \frac{b(b-1)(2b-1)}{6}
\]
with equality under the corresponding conditions.
The result follows from Corollary \ref{pairs} and the rearrangement inequality. The upper bound is $\sum_{k=1}^{b-1} k^2$ and the lower bound is $\sum_{k=1}^{b-1} k(b-k)$.

The second part is equivalent to showing that
\[
\frac{b(b-1)(b-2)}{6} \leq \tilde{S}_{(a_1,a_2;b)}(t) \leq \frac{b(b-1)(b-2)}{3}.
\]
for $1 \leq t \leq b-1$ with equality under the corresponding conditions.. Observe that as $k$ varies from $1$ to $b-1$, the terms $[a_2^{-1}(-t-k a_1)]$ vary over $\{0,1,\ldots,b-1\} \setminus \{x\}$ for some $1 \leq x \leq b-1$. By the rearrangement inequality,
\[
\sum_{k=1}^{b-1} (b-k)(k-1) \leq \sum_{k=1}^{b-1} k[a_2^{-1}(t-k a_1)] \leq \sum_{k=1}^{b-1}k(k-1).
\]
Simplifying, we obtain the result.
\proofend

We can easily translate the first statement to a statement about Dedekind sums, to recover a known bound. Indeed, $ S_{(a_1,a_2;b)}(0)=-s(a_1 a_2^{-1},b)+\frac{b-1}{4b}$, where $s$ is the Dedekind sum \cite{BR1}.

\begin{corollary}
Let $s(a,b)$ be the Dedekind sum. Then
\[
-\frac{(b-1)(b-2)}{12b} \leq s(a,b) \leq \frac{(b-1)(b-2)}{12b}. 
\]
The upper bound holds if and only if $a \equiv 1 \pmod{b}$ and the lower bound holds if and only if $a \equiv -1 \pmod{b}$.
\end{corollary}
 
 We see that $\{\tilde{S}_{(a_1,a_2;b)}(0)\in \R: (a_i,b)=1\}$ is a subset of
 \[
 \{(1,2,\ldots,b-1) \cdot (\sigma(1),\sigma(2),\ldots,\sigma(b-1)) \in \R: \sigma \in S_{b-1} \text{ is fixed-point free or is the identity}\}.
 \]
 To be more precise than that, we must know what permutation of $1+b\Z,2+b\Z,\ldots,(b-1)+b\Z$ is induced by multiplication by an integer relatively prime to $b$.

We now begin our pursuit of the second goal, to understand the location of maxima and minima of $2$-dimensional Fourier--Dedekind sums.
Define 
\[
\backslash \frac{x}{b}/:=\Big( \Big(\frac{x}{b} \Big) \Big)+\frac{1}{2} \delta_{\Z}(\frac{x}{b}).
\]
An equivalent definition is $\backslash \frac{x}{b}/:=\frac{x}{b}-\frac{1}{2}$, for $1 \leq x \leq b$ and extended periodically. This function arises naturally, e.g., $S_{(a;b)}(t)=\backslash \frac{a^{-1} t}{b}/-\frac{1}{2b}$, and will be convenient to work with in this section. We also mention in passing that Reciprocity for $S_{(a;b)}(t)$ is equivalent to Reciprocity of $\backslash \frac{a^{-1} t}{b}/$: 
\[
\backslash \frac{a^{-1} t}{b}/+\backslash \frac{b^{-1} t}{a}/=\frac{t}{ab}
\]
for $1 \leq t \leq a+b-1$. By a change of variables, it will suffice to understand the maxima and minima of $S_{(a,1;b)}(t)$.

\begin{theorem} \label{Optimization} For $a,b,t \in \N$,
\[
\argmaxl_{1 \leq t \leq b} S_{(a,1;b)}(t) \subset [\frac{b+1}{2},\frac{b+1}{2}+a].
\]
and
\[
\argminl_{1 \leq t \leq b} S_{(a,1;b)}(t) \subset [1,\min\{a,\frac{b+1}{2}\}].
\]
\end{theorem}

\proof
By Theorem \ref{functionalEqConv}, if $M$ is a maximizer, then for all $k \in \N$,
\[
S_{(a,1;b)}(M)-S_{(a,1;b)}(M+k a)=\sum_{j=0}^{k-1} S_{(1;b)}(M+ja) \geq 0.
\]
Recall that 
\[
S_{(1;b)}(t)= \backslash \frac{t}{b}/-\frac{1}{2b}.
\]
For $k=1$, we must have
\[
\backslash \frac{M}{b}/-\frac{1}{2b} \geq 0.
\]
\[
\iff M \geq \frac{b+1}{2}.
\]
Similarly, $S_{(a,1;b)}(M-a)-S_{(a,1;b)}(M)= S_{(1;b)}(M-a) \leq 0$. Consequently,
\[
\backslash \frac{M-a}{b}/-\frac{1}{2b} \leq 0.
\]
We show that $M \geq a+1$. Indeed, assume otherwise. Then we must have
\[
\frac{M-a}{b}-\frac{1}{2}+1-\frac{1}{2b}=\frac{M-a+\frac{b-1}{2}}{b} \leq 0.
\]
But by the case $k=1$, $M \geq \frac{b+1}{2}$. Therefore $M \geq a+1$. Using $\backslash \frac{M-a}{b}/-\frac{1}{2b} \leq 0$, we see that $M \leq \frac{b+1}{2}+a$.

Now let $1 \leq m \leq b$ be a minimizer. We have
\[
S_{(a,1;b)}(m)-S_{(a,1;b)}(m+a)=S_{(1;b)}(m) \leq 0.
\]
Consequently,
\[
\backslash \frac{m}{b}/-\frac{1}{2b}=\frac{m-\frac{b+1}{2}}{b} \leq 0.
\]
We also have
\[
S_{(a,1;b)}(m-a)-S_{(a,1;b)}(m)= S_{(1;b)}(m-a) \geq 0
\]
so that
\[
\backslash \frac{m-a}{b}/-\frac{1}{2b} \geq 0.
\]
Assume by contradiction that $m \geq a+1$. Then
\[
\backslash \frac{m-a}{b}/-\frac{1}{2b}=\frac{m-a-\frac{b+1}{2}}{b} \geq 0,
\]
a contradiction.
\proofend

\begin{figure}[H]
\centering
\includegraphics[width=5in]{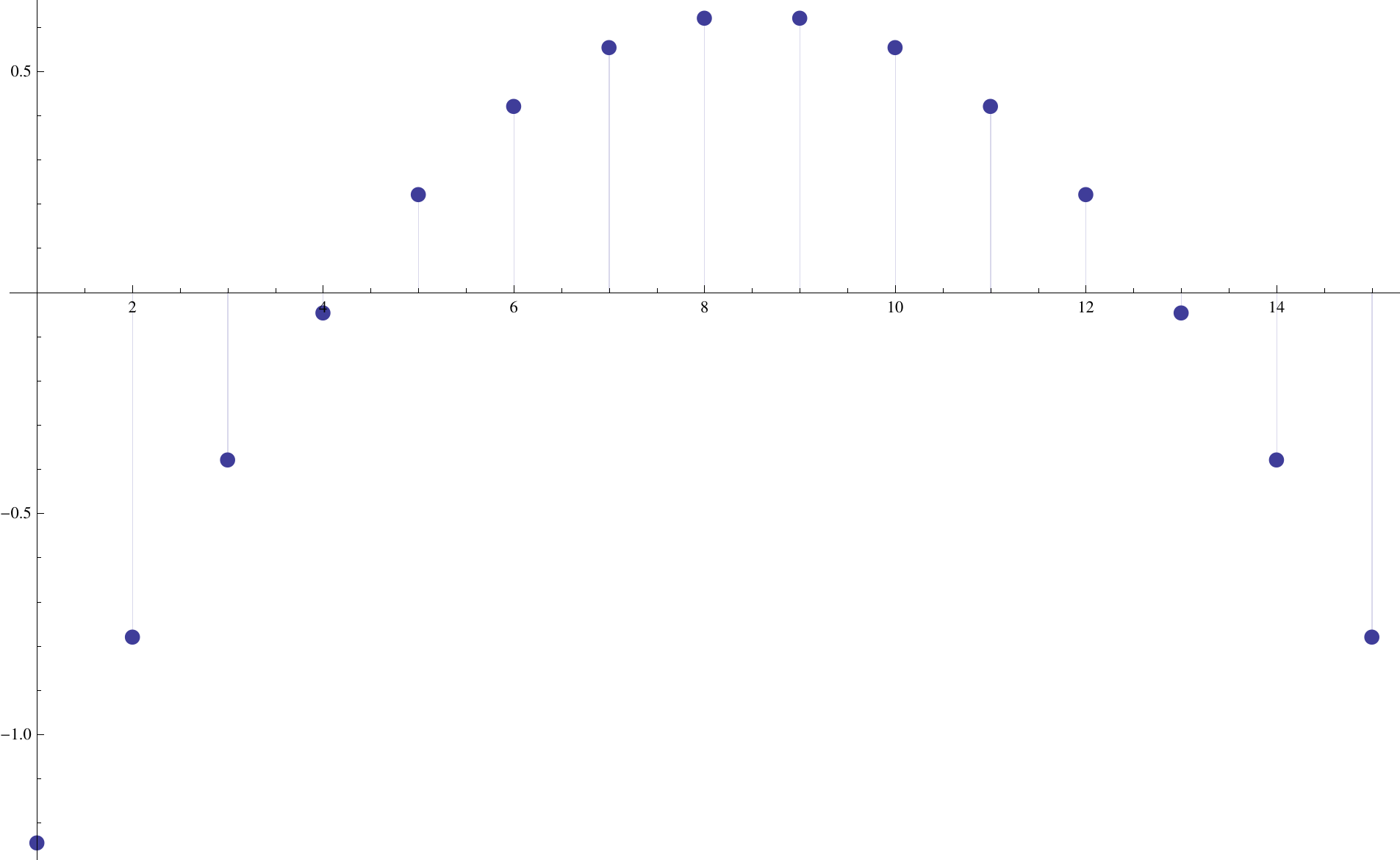}
\caption{\label{figargmax1} Plot of $S_{(1,1;15)}(t)$. Theorem \ref{Optimization} shows that the argmax is in the interval $[8,9]$ and the argmin is at $1$.}
\end{figure}

\begin{figure}[H]
\centering
\includegraphics[width=5in]{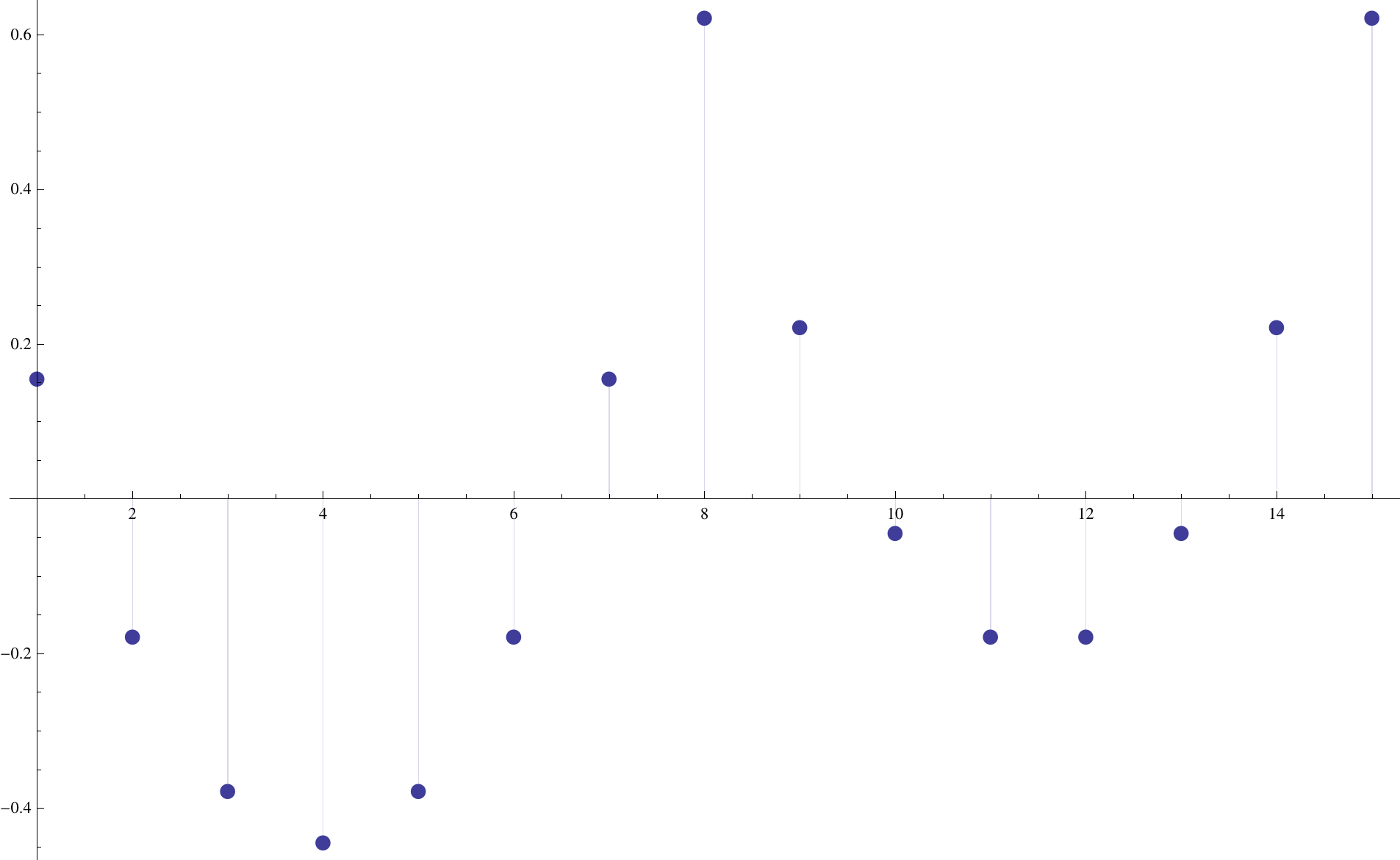}
\caption{Plot of $S_{(7,1;15)}(t)$. Theorem \ref{Optimization} shows that the argmax is in the interval $[8,15]$ and the argmin is in the interval $[1,7]$.}
\end{figure}

Theorem \ref{functionalEqConv} gives us not only information about the location of maxima and minima, but also about the concavity of the graph. Indeed,
\[
 (I-T^{a_{1}})(I-T^{a_{2}}) S_{(a_1,a_2;b)}=S_b=\delta_{\Z}(\frac{t}{b})-\frac{1}{b},
\]
by Lemma \ref{geometricSeries}. Consequently, 
\begin{proposition} \label{concavity} Let $a_1,a_2,b \in \N$ with $b$ relatively prime to $a_1,a_2$.
For $1 \leq t \leq b-1$,
\[
\frac{S_{(a_1,a_2;b)}(t)+S_{(a_1,a_2;b)}(t+a_1+a_2)}{2} > \frac{S_{(a_1,a_2;b)}(t+a_1)+S_{(a_1,a_2;b)}(t+a_2)}{2}
\]
and for $t=0$,
\[
\frac{S_{(a_1,a_2;b)}(t)+S_{(a_1,a_2;b)}(t+a_1+a_2)}{2} < \frac{S_{(a_1,a_2;b)}(t+a_1)+S_{(a_1,a_2;b)}(t+a_2)}{2}.
\]
\end{proposition}

\begin{figure}[H]
\centering
\includegraphics[width=5in]{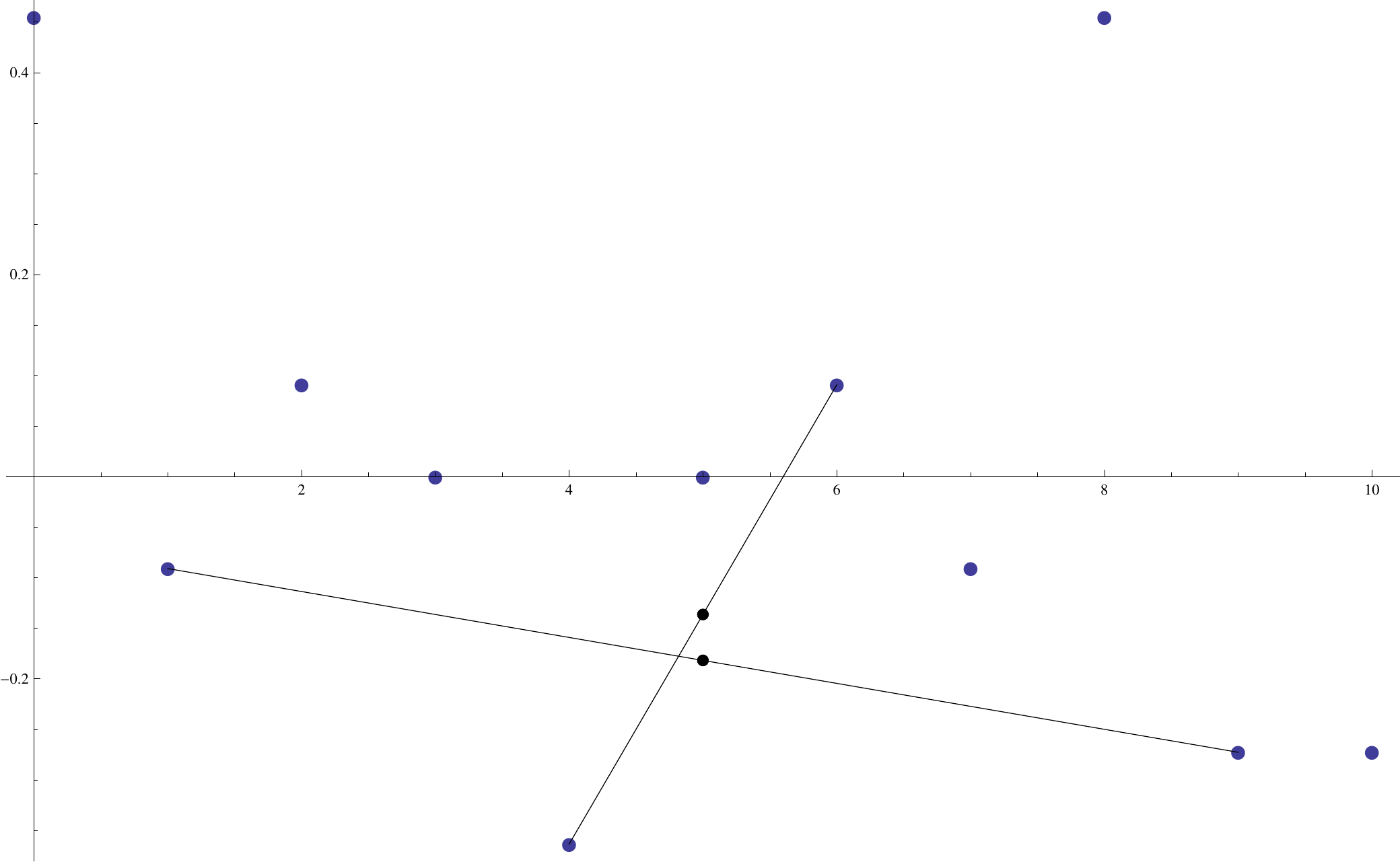}
\caption{\label{Concave} A plot of $S_{(3,5;11)}(t)$. Proposition \ref{concavity} implies that the average value of $S_{(3,5;11)}(4)$ and $S_{(3,5;11)}(6)$  is greater than the average value of $S_{(3,5;11)}(1)$ and $S_{(3,5;11)}(9)$.}
\end{figure}

Consequently, in the sense described in Proposition \ref{concavity}, the function $S_{(a_1,a_2;b)}(t)$ is mostly ``concave''. For an illustration, see Figures \ref{Concave} and \ref{figargmax1}.

We now seek to find good bounds on $R_{(a,1,b)}(t)$.  We note that by evaluating the lattice point enumerator in equation (\ref{L(t)eq}) at $t=0$, one can show that
\[
R_{(a,1,b)}(0)=1-\frac{1}{4}(1+\frac{1}{a}+\frac{1}{b})-\frac{1}{12}(\frac{a}{b}+\frac{b}{a}+\frac{1}{ab}).
\]

\begin{lemma} \label{kaplusb}
Let $k \in \N$. Then 
\[
R_{\{a,b\}}(t+k(a+b))=R_{\{a,b\}}(t)-[\sum_{j=0}^{k-1} \delta_{\Z}(\frac{t+jb}{a})+\delta_{\Z}(\frac{t+ja}{b})]+\frac{k}{a}+\frac{k}{b}.
\]
\end{lemma}

\proof
By Lemma \ref{RPIE},
\[
R_{\{a,b\}}(t+k(a+b))=R_{\{a,b\}}(t+(k-1)(a+b))-R_{a}(t+(k-1)(a+b))-R_{b}(t+(k-1)(a+b))=
\]
\[
=R_{\{a,b\}}(t+(k-1)(a+b))-R_{a}(t+(k-1)b)-R_{b}(t+(k-1)a)
\]
\[
=R_{\{a,b\}}(t+(k-1)(a+b))-\delta_{\Z}(\frac{t+(k-1)b}{a})-\delta_{\Z}(\frac{t+(k-1)a}{b})-\frac{1}{a}-\frac{1}{b}=
\]
\[
=R_{\{a,b\}}(t)-[\sum_{j=0}^{k-1} \delta_{\Z}(\frac{t+jb}{a})+\delta_{\Z}(\frac{t+ja}{b})]+\frac{k}{a}+\frac{k}{b}.
\]
\proofend

\begin{lemma} \label{deltaz} For all $t,m,a,b \in \N$,
\[
\sum_{j=0}^{t-1} \delta_{\Z}(\frac{j+ma}{b})=\lfloor \frac{ma+t-1}{b} \rfloor-\lfloor \frac{ma-1}{b} \rfloor.
\]
\end{lemma}

\proof
In general,
\[
\sum_{j=0}^{t-1} \delta_{\Z}(\frac{j}{b})=\lfloor \frac{t-1}{b} \rfloor+1.
\]
Subtracting yields the result.
\proofend

Let $k \in \N$. Applying Theorem \ref{functionalEqConv}, we see that
\[
R_{\{a,1,b\}}(k(a+b))-R_{\{a,1,b\}}(t+k(a+b))=\sum_{j=0}^{t-1}R_{\{a,b\}}(k(a+b)+j).
\]
By Proposition \ref{kaplusb}, this expression is equal to
\[
\sum_{j=0}^{t-1} R_{\{a,b\}}(j)+k(\frac{1}{a}+\frac{1}{b})-\sum_{m=0}^{k-1}[ \delta_{\Z}(\frac{j+mb}{a})+\delta_{\Z}(\frac{j+ma}{b})]= 
\]
\[
=R_{\{a,1,b\}}(0)-R_{\{a,1,b\}}(t)+tk(\frac{1}{a}+\frac{1}{b})-\sum_{j=0}^{t-1} \sum_{m=0}^{k-1} [\delta_{\Z}(\frac{j+mb}{a})+\delta_{\Z}(\frac{j+ma}{b})].
\]
Interchanging the order of summation and using Lemma \ref{deltaz}, we see that

\[
R_{\{a,1,b\}}(k(a+b))-R_{\{a,1,b\}}(t+k(a+b))= 
\]
\begin{align} \label{eq:diff}
R_{\{a,1,b\}}(0)-R_{\{a,1,b\}}(t)+tk(\frac{1}{a}+\frac{1}{b})-\sum_{m=0}^{k-1} \lfloor \frac{ma+t-1}{b} \rfloor-\lfloor \frac{ma-1}{b} \rfloor+\lfloor \frac{mb+t-1}{a} \rfloor-\lfloor \frac{mb-1}{a} \rfloor.
\end{align}

\begin{theorem} \label{boundstaplusb} For every $t \in \Z$,
\[
|R_{\{a,1,b\}}(t+a+b)-R_{\{a,1,b\}}(t)| \leq 1-\frac{1}{2}(\frac{1}{a}+\frac{1}{b}).
\]
\end{theorem}

\proof
Setting $k=1$ in equation \ref{eq:diff} gives

\[
R_{\{a,1,b\}}(a+b)-R_{\{a,1,b\}}(t+a+b)=
\]
\[
=R_{\{a,1,b\}}(0)-R_{\{a,1,b\}}(t)+t(\frac{1}{a}+\frac{1}{b})-\lfloor \frac{t-1}{b} \rfloor+\lfloor \frac{-1}{b} \rfloor-\lfloor \frac{t-1}{a} \rfloor+\lfloor \frac{-1}{a} \rfloor
\]
\[
=R_{\{a,1,b\}}(0)-R_{\{a,1,b\}}(t)+t(\frac{1}{a}+\frac{1}{b})-\lfloor \frac{t-1}{b} \rfloor-\lfloor \frac{t-1}{a} \rfloor-2.
\]

We show that for $t \in \Z$, 
\[
t(\frac{1}{a}+\frac{1}{b})-\lfloor \frac{t-1}{b} \rfloor-\lfloor \frac{t-1}{a} \rfloor-2 \leq 0.
\]
Multiplying both sides by $ab$ and rearranging, the statement is equivalent to
\[
2 ab \geq a((t-1)-b\lfloor \frac{t-1}{b} \rfloor)+b((t-1)+a\lfloor \frac{t-1}{a} \rfloor)+a+b.
\]
Recall that $x-y \lfloor \frac{x}{y} \rfloor = x \bmod y$. Since $t \in \Z$, 
\[
(t-1) \bmod b \leq b-1
\]
and
\[
(t-1) \bmod a \leq a-1.
\]
It follows that
\[
a((t-1)-b\lfloor \frac{t-1}{b} \rfloor)+b((t-1)+a\lfloor \frac{t-1}{a} \rfloor)+a+b \leq a(b-1)+b(a-1)+a+b=2ab.
\]
The upper bound is straightforward. Simplifying and using Reciprocity yields the result.
\proofend

\begin{figure}[H]
\centering
\includegraphics[width=5in]{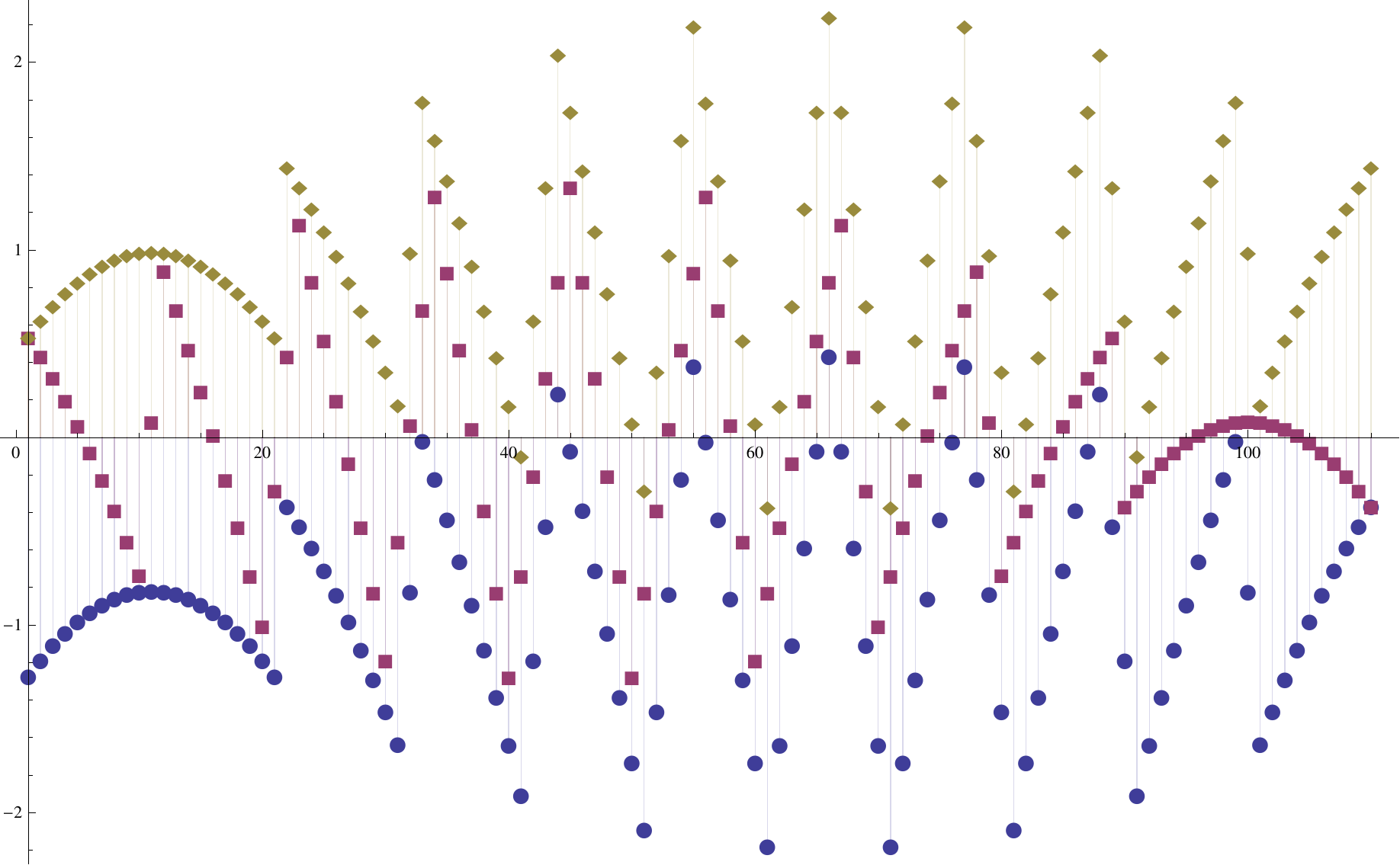}
\caption{Theorem \ref{boundstaplusb}. Here $a=11$ and $b=10$. We plot $R_{\{a,1,b\}}(t+a+b)$, the lower bound $R_{\{a,1,b\}}(t)-1+\frac{1}{2}(\frac{1}{a}+\frac{1}{b})$ and the upper bound $R_{\{a,1,b\}}(t)+1-\frac{1}{2}(\frac{1}{a}+\frac{1}{b})$ over $t=1,2,\ldots,a b$.}
\end{figure}

\begin{remark}
The inequality $t(\frac{1}{a}+\frac{1}{b})-\lfloor \frac{t-1}{b} \rfloor-\lfloor \frac{t-1}{a} \rfloor-2 \leq 0$  in the proof of Theorem \ref{boundstaplusb} does not necessarily hold for a general real value of $t$.
\end{remark}

\begin{corollary} \label{boundsRecip}
For $t=1,2,\ldots,a+b$,
\[
|S_{(a,1;b)}(t+a+b)+S_{(b,1;a)}(t+a+b)+\frac{t^2}{2 a b}-\frac{t}{2}(\frac{1}{a}+\frac{1}{b}+\frac{1}{ab})+\frac{1}{12}(\frac{1}{a}+\frac{1}{b}+3+\frac{1}{ab}+\frac{a}{b}+\frac{b}{a})| \leq 1-\frac{1}{2}(\frac{1}{a}+\frac{1}{b})
\]
\end{corollary}

\proof
By Reciprocity, $R_{\{a,1,b\}}(t)=\frac{t^2}{2 a b}-\frac{t}{2}(\frac{1}{a}+\frac{1}{b}+\frac{1}{ab})+\frac{1}{12}(\frac{1}{a}+\frac{1}{b}+3+\frac{1}{ab}+\frac{a}{b}+\frac{b}{a})$.
\proofend

\begin{figure}[H]
\centering
\includegraphics[width=5in]{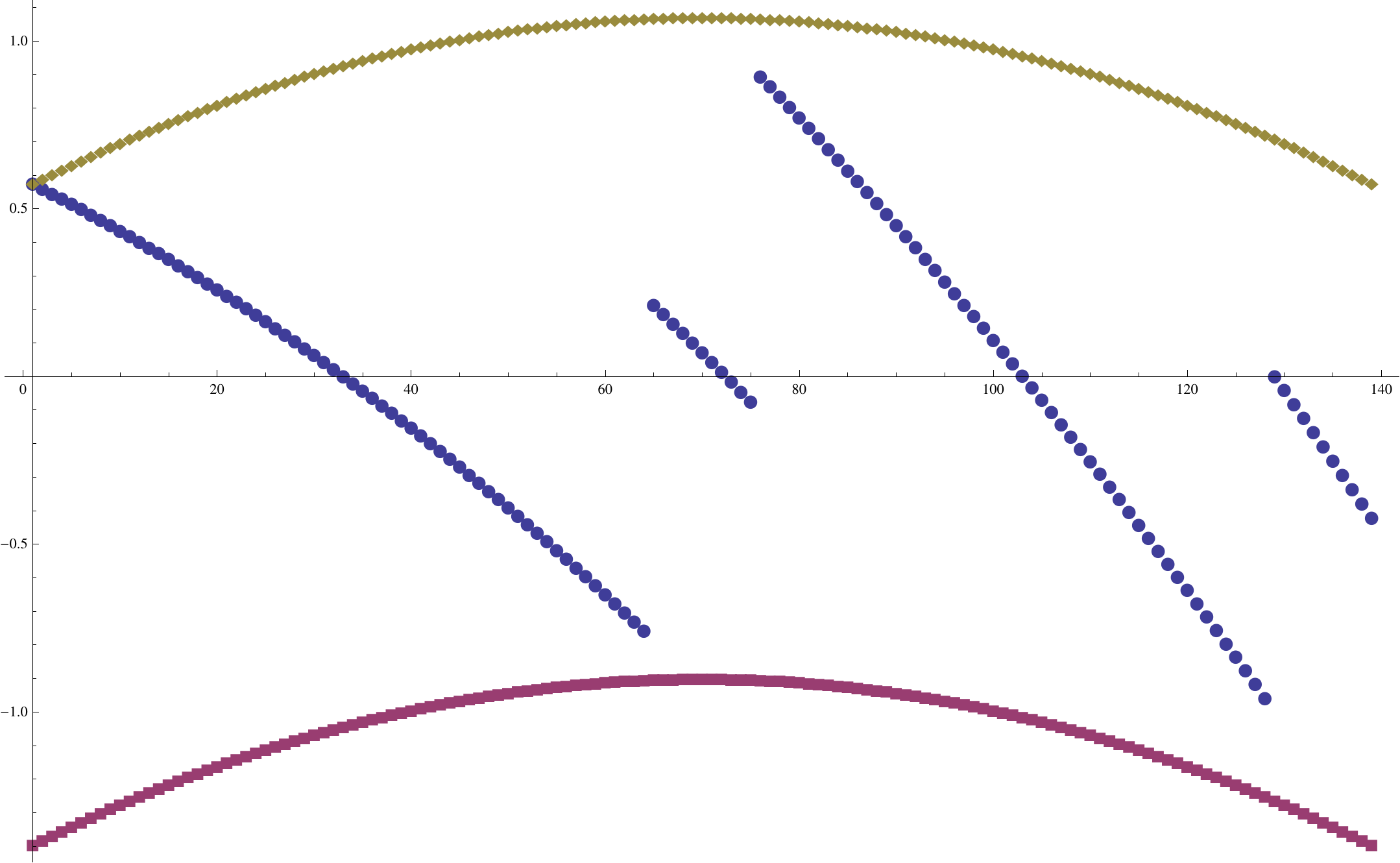}
\caption{An example of the bounds in Corollary \ref{boundsRecip} for $a=64$ and $b=75$. Here we plot $R_{\{a,1,b\}}(t+a+b)$ and its lower and upper bounds over $t=1,2,\ldots,a+b$.}
\end{figure}

\bigskip
{\bf Acknowledgments}. It is a pleasure to acknowledge interesting discussions and valuable feedback from M. Beck, D. Cristofaro-Gardiner, K. Ribet, S. Robins, J. Sondow and X. Yuan.
This work originated during the 2013 UC Berkeley Geometry, Topology, and Operator Algebras RTG Summer Research Program for Undergraduates. This program was funded by the National Science
Foundation.

This material is based upon work supported by the National Science Foundation Graduate Research Fellowship under Grant No. DGE 1106400. Any opinion, findings, and conclusions or recommendations expressed in this material are those of the authors(s) and do not necessarily reflect the views of the National Science Foundation.


\begin{thebibliography}{99}
\bibitem{BDR1} M. Beck, R. Diaz, and S. Robins. \textit{The Frobenius problem, rational
polytopes, and Fourier-Dedekind sums.} J. Number Theory, 96(1):1-21, 2002.

\bibitem{B1} M. Beck, I. Gessel and T. Komatsu, \textit{The polynomial part of a restricted partition function related to the Frobenius problem,} Electron. J. Combin. \textbf{8}, No. 1 (2001), N 7.

\bibitem{BR1} M. Beck and S. Robins, \textit{Computing the continuous discretely: Integer-point enumeration in
polyhedra.} Undergraduate Texts in Mathematics. Springer, New York, 2007.

\bibitem{HZ1} F. Hirzebruch, D. Zagier, \textit{The Atiyah-Singer theorem and elementary number theory,}
Publish or Perish, Boston (1974).

\bibitem{Hutch1} M. Hutchings, \textit{Recent progress on symplectic embedding problems in four
dimensions,} Proc. Natl. Acad. Sci. USA \textbf{108} (2011), 8093-8099.

\bibitem{K1} D. Knuth, \textit{The art or computer programming, vol. 2,} Addison-Wesley, Reading, Mass.,
(1981).

\bibitem{MD1} D. McDuff, \textit{The Hofer conjecture on embedding symplectic ellipsoids,} J.
Diff. Geom. 88 (2011), 519-532.





\end{thebibliography}
\end{document}